\documentclass{article}

\usepackage{amssymb}
\usepackage{amsmath}
\usepackage{amsthm}
\usepackage{tikz-cd}
\usepackage{enumerate}

\numberwithin{equation}{section}
\renewcommand{\labelenumi}{(\roman{enumi})}

\newtheoremstyle{slplain}% name
  {\topsep}%Space above
  {\topsep}%Space below
  {\slshape}%Body font
  {0pt}%Indent amount
  {\bfseries}% Theorem head font
  {.}%Punctuation after theorem head
  {0.5em}%Space after theorem head 2
  {}%Theorem head spec (can be left empty, meaning ‘normal’)

\theoremstyle{slplain}
  \newtheorem{THM}{Theorem}[section]
  \newtheorem{LEM}[THM]{Lemma}
  \newtheorem{COR}[THM]{Corollary}
  \newtheorem{ASSUMP}[THM]{Assumption}
  \newtheorem*{OPEN}{Open problem}
  \newtheorem*{THMNONUM}{Theorem}
  \newtheorem*{CORNONUM}{Corollary}

\theoremstyle{definition}
  \newtheorem{DEF}[THM]{Definition}
  \newtheorem{EX}{Example}[section]
  
  \newtheorem{CONSTR}[THM]{Construction}

\newcommand\AAA{\mathbf{A}}
\newcommand\BB{\mathbf{B}}
\newcommand\CC{\mathbf{C}}
\newcommand\DD{\mathbf{D}}
\newcommand\EM{\mathbf{EM}}
\newcommand\KK{\mathbf{K}}
\newcommand\NN{\mathbb{N}}
\newcommand\QQ{\mathbb{Q}}
\newcommand\VV{\mathbf{V}}
\newcommand\WW{\mathbf{W}}

\newcommand\calA{\mathcal{A}}
\newcommand\calB{\mathcal{B}}
\newcommand\calC{\mathcal{C}}

\newcommand\calF{\mathcal{F}}
\newcommand\calH{\mathcal{H}}

\newcommand\calS{\mathcal{S}}
\newcommand\calX{\mathcal{X}}
\newcommand\calY{\mathcal{Y}}

\newcommand\0{\varnothing}
\renewcommand\epsilon{\varepsilon}
\renewcommand\phi{\varphi}
\let\le\leqslant
\let\ge\geqslant
\newcommand{\lt}[1]{\mathrel{<^{#1}}}

\newcommand{\lex}[1]{\mathrel{<^{#1}_{\mathit{lex}}}}

\newcommand{\restr}[2]{\mbox{$#1$}\mbox{$\upharpoonright$}_{#2}}

\newcommand\Ob{\mathrm{Ob}}

\newcommand\ID{\mathrm{ID}}
\newcommand\id{\mathrm{id}}

\newcommand{\Boxed}[1]{\hbox{$#1$}}
\newcommand{\Borel}[1]{\overset{\flat}{#1}}
\newcommand{\flatT}{T^{\flat}}
\newcommand{\flatdualT}{T^{\flat\partial}}

\newcommand{\Set}{\mathbf{Set}}

\newcommand{\Top}{\mathbf{Top}}
\newcommand{\ChEmb}{\mathbf{Ch}_{\emb}}
\newcommand{\ChEmbFin}{\mathbf{Ch}_{\emb}^{\fin}}
\newcommand{\ChRs}{\mathbf{Wch}_{\rs}}
\newcommand{\ChRsFin}{\mathbf{Wch}_{\rs}^{\fin}}
\newcommand{\WalgRe}{\mathbf{Walg}_{\re}}
\newcommand{\WalgReFin}{\mathbf{Walg}_{\re}^{\fin}}

\newcommand{\Alg}{\mathbf{Alg}}

\newcommand{\AlgEpi}{\mathbf{Alg}_{\epi}}

\newcommand{\VEpi}{\mathbf{V}_{\epi}}

\newcommand{\VEpiFin}{\mathbf{V}_{\epi}^{\fin}}

\newcommand\eval{\mathrm{eval}}
\newcommand\op{\mathit{op}}

\newcommand\fin{\mathit{fin}}
\newcommand\re{\mathit{re}}
\newcommand\rs{\mathit{rs}}
\newcommand\emb{\mathit{emb}}
\newcommand\epi{\mathit{epi}}
\newcommand\Aut{\mathrm{Aut}}
\newcommand\iso{\mathrm{iso}}

\newcommand\Fraisse{Fra\"\i ss\'e}

\title{Dual Ramsey properties for classes of algebras}
\author{Dragan Ma\v sulovi\'c\\
        Department of Mathematics and Informatics\\
        Faculty of Sciences, University of Novi Sad, Serbia\\
        email: dragan.masulovic@dmi.uns.ac.rs}

\begin{document}
\maketitle

\begin{abstract}
  Almost any reasonable class of finite relational structures has the Ramsey property or a precompact Ramsey expansion.
  In contrast to that, the list of classes of finite algebras with the precompact Ramsey expansion is surprisingly short.
  In this paper we show that any nontrivial variety (that is, equationally defined class of algebras)
  enjoys various \emph{dual} Ramsey properties. We develop a completely new set of strategies that rely on
  the fact that left adjoints preserve the dual Ramsey property,
  and then treat classes of algebras as Eilenberg-Moore categories for a monad.
  We show that finite algebras in any nontrivial variety have finite dual small Ramsey degrees,
  and that every finite algebra has finite dual big Ramsey degree in the free algebra on countably many free generators.
  As usual, these come as consequences of ordered versions of the statements.

  \bigskip


  \bigskip

  \noindent \textbf{AMS Subj.\ Classification (2020):} 05C55; 18C20
\end{abstract}

\section{Introduction}

Almost any reasonable class of finite relational structures has the Ramsey property or a precompact Ramsey expansion
(see~\cite{vanthe-more} for definition), usually by an appropriate choice of linear orders. For example,
finite graphs expanded with arbitrary linear orders have the Ramsey property~\cite{Nesetril-Rodl,Nesetril-Rodl-2};
the same holds for finite hypergraphs~\cite{Nesetril-Rodl,Nesetril-Rodl-2} and finite metric spaces~\cite{Nesetril-metric};
finite posets expanded with linear orders that extend the partial order have the Ramsey property~\cite{PTW,Fouche-1997};
the same holds for multiposets --- structures with several partial orders forming a partial order~\cite{draganic-masul};
finite equivalence relations with linear orders where equivalence classes are convex have the Ramsey property~\cite{KPT};
the same holds for finite ultrametric spaces with linear orders where balls are convex~\cite{van-the-metric-spaces,van-the-memoirs};
and the list goes on and on. One of the most prominent general results in this direction is the Ne\v set\v ril-R\"odl Theorem:

\begin{THM}[The Ne\v set\v ril-R\"odl Theorem \cite{Nesetril-Rodl,Nesetril-Rodl-2}]\label{rpemklei.thm.N-R}
  Consider a finite relational language with a distinguished binary relational symbol $<$ which
  is always interpreted as a linear order.
  Let $F$ be a set of finite ordered irreducible\footnote{An ordered relational structure is \emph{irreducible} if every pair of
  distinct elements of the structure appears in some tuple of some relation in the language distinct from the distinguished
  binary relation.} relational structures over the language and let $\KK$ be the class of all the finite ordered relational
  structures that embed no finite structure from $F$. Then $\KK$ has the Ramsey property.
\end{THM}

The Ramsey property imposes severe restrictions on the classes of objects enjoying the property.
One of those restrictions is the requirement that all the objects in such classes be rigid.
This was observed for structures in~\cite{Nesetril,zucker1} and generalized to categorical setting in~\cite{masul-scow}
(see Theorem~\ref{rpemklei---rigid} below). Hence, natural classes of structures such as finite graphs and finite posets
usually do not have the Ramsey property. Nevertheless, many of these classes enjoy
the weaker property of having finite (small) Ramsey degrees. Zucker recently proved
(see~\cite[Theorem~8.14]{zucker1}) that a class of finite structures has small Ramsey degrees
if and only if it has a precompact Ramsey expansion. The phenomenon of having a precompact Ramsey expansion
(or, equivalently, finite small Ramsey degrees) is so ubiquitous that it is generally believed that
every reasonable class of finite relational structures has a precompact Ramsey expansion, though not necessarily by linear
orders (see \cite{vanthe-more}).

In contrast to that, the list of classes of finite algebras with a precompact Ramsey expansion is surprisingly short.
The Ramsey property has been established for the following classes of finite algebras:
finite boolean algebras~\cite{graham-rothschild-DRT}; finite vector spaces over a
fixed finite field~\cite{graham-leeb-rothschild}; finite boolean lattices~\cite{promel-85};
finite unary algebras over a finite language~\cite{sokic-unary-functions};
finite $G$-sets for a finite group $G$~\cite{sokic-unary-functions};
and finite semilattices~\cite{sokic-semilattices}. One of the most prominent general results in this direction,
the Theorem of Evans, Hubi\v cka and Ne\v set\v ril, treats the operations in the language as special relations requiring, thus,
that we allow for partial operations.

\begin{THM}[The Evans-Hubi\v cka-Ne\v set\v ril Theorem~\cite{evans-hubicka-nesetril}]
  Every free amalgamation class of ordered first-order structures over the same first-order language
  has the Ramsey property (where functional symbols in the language are interpreted as partial operations).
\end{THM}

A closer inspection reveals that there are not many natural (equationally defined) classes of algebras
among the classes of algebras where we have identified the Ramsey property.
%And this is not for lack of trying!
%Just one quick look at~\cite{evans-hubicka-nesetril} or~\cite{hubicka-nesetril-all-those}
%suffices to see that extremely sophisticated machinery and almost superhuman endurance
%was needed to achieve these results.
For example, it was shown in~\cite{kechris-sokic} that no expansion of the class of finite
distributive lattices by linear orders satisfies the Ramsey property (although there is an expansion using ternary relations,
see~\cite{kechris-sokic}).
Motivated by this result in~\cite{masul-WTC} we show that for an arbitrary nontrivial locally finite variety $\VV$
of lattices distinct from the variety of all lattices and the variety of distributive lattices,
no reasonable expansion of $\VV^\fin$ ($=$ the class of all the finite lattices in $\VV$) has the Ramsey property.
However, if we consider lattices as partially ordered sets (and thus switch from
the lattices as algebras to their relational alter ego)
we show in \cite{masul-WTC} that \emph{every} variety of lattices
gives rise to a class of linearly ordered posets having both the Ramsey property and
the ordering property (see~\cite{Nesetril} for definition).
It seems that structural Ramsey theory is at odds with natural (equationally defined) classes of algebras.

This point of view is further supported by the following question which, despite a long
history~\cite{deuber-rothschild,oates-williams,voigt-fin-ab-grps}, is still open:

\begin{OPEN}
  Is it true that the class of finite groups has a precompact Ramsey expansion?
\end{OPEN}

In this paper we show that nontrivial varieties of algebras enjoy various \emph{dual} Ramsey properties.
The search for dual Ramsey statements has been an important research direction in the past 50 years not only
because dual Ramsey results are relatively rare in comparison to the vast number of ``direct'' Ramsey results, but also
because they require intricate proof strategies and are usually more powerful than their ``direct'' analogues.
It turns out that classes of algebras are a gold mine of dual Ramsey results, and one of the first papers to point in this
direction is~\cite{solecki-2010}.

In contrast to most mainstream Ramsey-related results where the context is that of model theory (see for example~\cite{solecki-2010}),
all our results are spelled out using the categorical reinterpretation of the Ramsey property as proposed in~\cite{masul-scow}.
Actually, it was Leeb who pointed out already in early 1970's that the use of category theory can be quite helpful
both in the formulation and in the proofs of results pertaining to structural Ramsey theory~\cite{leeb-cat}.
We strongly believe that this is even more the case when dealing with the dual Ramsey property.

In order to prove various dual Ramsey statements for classes of algebras
we develop a completely new set of strategies that rely on
the fact that right adjoints preserve the Ramsey property while left adjoints preserve the dual Ramsey property.
We then consider varieties (that is, equationally defined classes) of algebras
as Eilenberg-Moore categories for a monad and show the following:

\begin{THMNONUM}[see Theorem~\ref{rpemklei.thm.ALG-MAIN} below and also~\cite{solecki-2010}]
  Let $\VV$ be a nontrivial variety of algebras over the same algebraic language and
  $\KK$ the class of all ordered finite algebras from $\VV$. Then $\KK$ has the dual Ramsey property
  with respect to rigid epimorphisms (that is, epimorphisms that are at the same time rigid surjections).
\end{THMNONUM}

A proof of the same result that relies on model-theoretic strategies can be found in~\cite{solecki-2010}.
Using further categorical properties of small Ramsey degrees
we can then get rid of linear orders to derive the following statement:

\begin{THMNONUM}[see Corollary~\ref{rpemklei.cor.ALG-MAIN} below]
  For every nontrivial variety $\VV$ of algebras over the same algebraic language, finite algebras in
  $\VV$ have finite dual small Ramsey degrees with respect to epimorphisms.
\end{THMNONUM}

In comparison to the Ne\v set\v ril-R\"odl Theorem (Theorem~\ref{rpemklei.thm.N-R})
which tells us that every structured class of finite ordered relational structures has
the Ramsey property, the above result can be thought as a ``result through the looking-glass'': every
structured class of finite ordered algebras has the dual Ramsey property.
In particular, by specializing to groups the following conclusion is straightforward:

\begin{CORNONUM}[see Corollary~\ref{rpemklei.cor.groups} below]
  The class of all finite ordered groups has the dual Ramsey property. Moreover, for every nontrivial variety of
  groups (abelian groups, for example), the class of ordered finite groups from the variety has the dual Ramsey property
  with respect to rigid epimorphisms.

  Consequently, for every nontrivial variety $\VV$ of groups, finite groups in
  $\VV$ have finite dual small Ramsey degrees with respect to epimorphisms.
\end{CORNONUM}

As the language of small Ramsey degrees enables us to talk about the Ramsey property
in the context of finite structures that are not rigid, the language of big Ramsey degrees makes it possible
to consider the Ramsey property of finite non-rigid structures with respect to an infinite universal structure.
The infinite version of Ramsey's Theorem \cite{ramsey} can be understood as the first result in this direction:
every finite chain has finite big Ramsey degree in $\omega$ (and that the degree is~1;
note that $\omega$ is universal for the class of all finite chains).

The study of big Ramsey degrees was explicitly suggested for the first time in~\cite{KPT}, although founding ideas
date back to the work of Galvin in late 1960's~\cite{galvin1,galvin2}. Big Ramsey degrees of finite chains in $\QQ$
were computed by Devlin in~\cite{devlin}. Sauer proved in \cite{Sauer-2006} that every finite graph has finite
big Ramsey degree in the Rado graph --- the \Fraisse\ limit of the class of all the finite graphs.
Nugyen Van Th\'e proved in~\cite{van-the-bigrd-umet} that for every nonempty finite set $S$ of non-negative reals,
every finite $S$-ultrametric space has finite big Ramsey degree in the \Fraisse\ limit of the class of all
the finite $S$-ultrametric spaces. Laflamme, Nugyen Van Th\'e and Sauer proved in~\cite{laf-vanthe-sauer-2010}
that every finite local order has finite big Ramsey degree in the dense local order $\calS(2)$.
A suite of remarkable results of Dobrinen \cite{dobrinen1,dobrinen2,dobrinen3} shows that 
every finite $K_n$-free graph has finite big Ramsey degree in the Henson graph $\calH_n$, and in some cases
the exact numbers can be produced. Finally, let us recall a result due to Zucker which is analogous to the
Ne\v set\v ril-R\"odl Theorem (Theorem~\ref{rpemklei.thm.N-R}):

\begin{THM}[Zucker's Theorem \cite{Zucker-NoteOnBigRD}]
  Fix a finite binary relational language. Let $F$ be a finite set of finite irreducible
  relational structures over the language and let $\KK$ be the class of all finite relational
  structures that embed no structure from $F$. Then
  every structure from $\KK$ has finite big Ramsey degree in the \Fraisse\ limit of $\KK$.
  (Note that $\KK$ is obviously an amalgamation class.)
\end{THM}

The final result of the paper is a ``looking-glass'' analogue of Zucker's Theorem.
At the very end of the paper we prove the following:

\begin{THMNONUM}[see Theorem~\ref{rpemklei.thm.main-V-order} and Corollary~\ref{rpemklei.cor.dual-big-rd-unordered} below]
  Let $\VV$ be a nontrivial variety of algebras over a countable algebraic language.
  Then every finite $\VV$ algebra has a finite big dual Ramsey degree in the free
  $\VV$ algebra on $\omega$ generators with respect to Borel colorings.
  Moreover, the big Ramsey degree of an algebra with $n$ elements does not exceed~$n \cdot n!$.
\end{THMNONUM}

The feeling that ``direct'' structural Ramsey theory is at odds with equationally defined classes of algebras
extends to considerations of big Ramsey degrees. For example,
we can consider the \Fraisse\ limit of finite boolean algebras (with embeddings as morphisms)
whose limit is the countable atomless boolean algebra. The following is an important open problem:

\begin{OPEN}
  Is it true that every finite boolean algebra has finite big Ramsey
  degree in the countable atomless boolean algebra?
\end{OPEN}

Although our main results strongly resemble their ``looking-glass'' analogues, the tools we need for their proofs are
of a completely different nature.
After fixing standard notions and notation in Section~\ref{rpemklei.sec.prelim}, we present our proof strategies in
Section~\ref{rpemklei.sec.RP}. Our starting point is the observation from~\cite{masul-scow} that right adjoints preserve the
Ramsey property while left adjoints preserve the dual Ramsey property. We then show that if $T$ is a monad on a
category with the dual Ramsey property, both the Kleisli category and the Eilenberg-Moore category for the monad have the
dual Ramsey property. Unfortunately, these two simple results are not very useful:
for the categorical treatment of the dual Ramsey property
it is essential to restrict the attention to categories where the morphisms are epi, and units for monads that
produce constructions we are interested in seldom enjoy the property.
Therefore, we relax the context by proving that
the dual Ramsey property carries over from a category to the category of \emph{weak Eilenberg-Moore algebras}
defined for functors with multiplication, which are straightforward weakenings of monads.
Finally, in Section~\ref{rpemklei.sec.drp-V} we prove the main results of the paper.
We show that for every algebraic language $\Omega$ and every nontrivial variety $\VV$ of $\Omega$-algebras
the class of finite ordered $\VV$ algebras taken with rigid epimorphisms (that is, epimorphisms of algebras that are at the
same time rigid surjections) has the dual Ramsey property. The unordered version then follows
immediately: for every algebraic language $\Omega$ and every nontrivial variety $\VV$ of $\Omega$-algebras,
finite $\VV$ algebras have finite small dual Ramsey degrees with respect to epimorphisms.
We then prove that for every countable algebraic language $\Omega$ and every nontrivial variety $\VV$ of $\Omega$-algebras
finite $\VV$ algebras have finite big Ramsey degrees in the free $\VV$ algebra
$\calF_\VV(\omega)$ on $\omega$ generators with respect to Borel colorings. As usual, this comes as a consequence of the
ordered version of the statement.

\section{Preliminaries}
\label{rpemklei.sec.prelim}

\paragraph{Chains and rigid surjections.}
  A \emph{chain} is a linearly ordered set $(A, \Boxed<)$.
  Finite or countably infinite chains will sometimes be denoted as
  $\{a_1 < \dots < a_n < \dots\}$. For example, $\omega = \{0 < 1 < 2 < \dots\}$.
  Every strict linear order $<$ induces the reflexive version $\le$ in the obvious way.
  We assume the Axiom of Choice so that every set can be well-ordered.

  Let $\alpha$ be an ordinal and let $(A_\xi, \Boxed<)_{\xi < \alpha}$, be a family of chains indexed by~$\alpha$.
  The \emph{ordinal sum} of the family $(A_\xi, \Boxed<)_{\xi < \alpha}$,
  is a new chain $\bigoplus_{\xi < \alpha} (A_\xi, \Boxed{<})$ constructed on a disjoint union
  $\bigcup_{\xi < \alpha} (\{\xi\} \times A_\xi)$ so that $(\xi, a) < (\eta, b)$ if $\xi < \eta$, or $\xi = \eta$ and $a < b$.
  The \emph{lexicographic product} of a \emph{finite} family $(A_i, \Boxed<)_{i < n}$
  is a new chain $\bigotimes_{i < n} (A_i, \Boxed{<})$ constructed on
  $\prod_{i < n} A_i$ as follows.
  For a pair of distinct elements $(a_0, \ldots, a_{n-1}) \ne (b_0, \ldots, b_{n-1})$ from the product
  let $\eta = \min \{i < n: a_i \ne b_i\}$. Then put $(a_0, \ldots, a_{n-1}) \lex{} (b_0, \ldots, b_{n-1})$
  if and only if $a_i < b_i$. Note that if all the chains $(A_i, \Boxed<)$, $i < n$, are well-ordered then
  $\bigotimes_{i < n} (A_i, \Boxed{<})$ is well-ordered.
  In particular, for every $n \in \NN$ and every well-ordered chain $(A, \Boxed{<})$ the chain $(A^n, \Boxed{\lex{}})$
  is also well-ordered.

  Let $(A, \Boxed{<})$ and $(B, \Boxed{<})$ be well-ordered chains.
  A surjective map $f : A \to B$ is a \emph{rigid surjection} from $(A, \Boxed{<})$ onto $(B, \Boxed{<})$ if 
  $b_1 < b_2$ implies $\min f^{-1}(b_1) < \min f^{-1}(b_2)$ for all $b_1, b_2 \in B$.

\paragraph{Algebras and varieties.}
  Let $\Omega$ be an algebraic language, that is, the set of constant and functional symbols.
  A \emph{$\Omega$-algebra} is a structure $(A, \Omega^A)$ where $\Omega^A = \{f^A : f \in \Omega\}$ is a set of operations
  on $A$ such that the arity of each operation $f^A$ coincides with
  the arity of the corresponding functional symbol $f \in \Omega$.
  For a class $\KK$ of $\Omega$-algebras let $\KK^\fin$ denote the class of all finite members of $\KK$.

  Recall that and algebra $\calF = (F, \Omega^F)$ is \emph{freely
  generated by a set of generators $X$} if $X \subseteq F$ and for every other algebra $\calA = (A, \Omega^A)$
  every set mapping $g : X \to A$ uniquely extends to a homomorphism $g^\# : \calF \to \calA$.

  Given an algebraic language $\Omega$ and a nonempty set of variables $X$ let $T(X)$ denote the set of all the $\Omega$-terms
  over the set of variables $X$. It is the carrier of the \emph{absolutely free algebra} $\calF(X) = (T(X), \Omega^{T(X)})$.
  In any algebra $\calA = (A, \Omega^A)$ each term $t \in T(X)$ in $n$ variables
  determines a function $t^A : A^n \to A$. Let $t_1, t_2 \in T(X)$ be terms in the same number of variables.
  An algebra $\calA$ satisfies an \emph{identity} $t_1 \approx t_2$, in symbols $\calA \models t_1 \approx t_2$,
  if $t_1^A = t_2^A$. A class of algebras $\KK$ satisfies the identity
  $t_1 \approx t_2$, in symbols $\KK \models t_1 \approx t_2$, if $\calA \models t_1 \approx t_2$ for all $\calA \in \KK$.

  A class $\VV$ of $\Omega$-algebras is a \emph{variety} if there is a set
  $\Sigma = \{t_1^i \approx t_2^i : i \in I\}$ of $\Omega$-identities such that an
  $\Omega$-algebra $\calA$ belongs to $\VV$ if and only if $\calA \models t_1^i \approx t_2^i$
  for all $i \in I$. A variety $\VV$ is \emph{nontrivial} if there exists an algebra $\calA = (A, \Omega^A) \in \VV$ such that $|A| \ge 2$.
  
  Given a variety $\VV$ of $\Omega$-algebras and a set of variables $X$ let
  $\Theta_\VV(X) = \{(t_1, t_2) \in T(X)^2 : \VV \models t_1 \approx t_2\}$. Clearly, $\Theta_\VV(X)$ is a congruence of
  the term algebra $T(X)$ and $T_\VV(X) = T(X) / \Theta_\VV(X)$ is the carrier of the \emph{free $\VV$ algebra $\calF_\VV(X)$ with
  free generators $X / \Theta_\VV(X)$}. Let $\nu_{\VV,X} : \calF(X) \to \calF_\VV(X)$ be the \emph{natural epimorphism}
  that takes $x$ to its equivalence class $x / \Theta_\VV(X)$.

\paragraph{Categories and functors.}
Let us quickly fix some basic category-theoretic notions and notation.
For a detailed account of category theory we refer the reader to~\cite{AHS}.

In order to specify a \emph{category} $\CC$ one has to specify
a class of objects $\Ob(\CC)$, a class of morphisms $\hom_\CC(A, B)$ for all $A, B \in \Ob(\CC)$,
the identity morphism $\id_A$ for all $A \in \Ob(\CC)$, and
the composition of morphisms~$\cdot$~so that
$\id_B \cdot f = f = f \cdot \id_A$ for all $f \in \hom_\CC(A, B)$, and
$(f \cdot g) \cdot h = f \cdot (g \cdot h)$ whenever the compositions are defined.
We write $A \overset\CC\longrightarrow B$ as a shorthand for $\hom_\CC(A, B) \ne \0$.

Given a category $\CC$, the \emph{opposite category} $\CC^\op$ is a category constructed from $\CC$ on the same class of objects
by formally reversing arrows and products. More precisely, for $A, B \in \Ob(\CC) = \Ob(\CC^\op)$ we have that
$\hom_{\CC^\op}(A, B) = \hom_{\CC}(B, A)$, and for $f \in \hom_{\CC^\op}(A, B)$ and $g \in \hom_{\CC^\op}(B, C)$
we have that $g \mathbin{\underset{\CC^\op}{\cdot}} f = f \mathbin{\underset{\CC}{\cdot}} g$.

A category $\CC$ is \emph{locally small} if $\hom_\CC(A, B)$ is a set for all $A, B \in \Ob(\CC)$.
Sets of the form $\hom_\CC(A, B)$ are then referred to as \emph{hom-sets}.
Hom-sets in $\CC$ will be denoted by $\hom_\CC(A, B)$, or simply $\hom(A, B)$ when $\CC$ is clear from the context.
All the categories in this paper are locally small. We shall explicitly state this assumption in the
formulation of main results, but may omit the explicit statement of this fact in the formulation
of auxiliary statements.

A morphism $f$ is: \emph{mono} or \emph{left cancellable} if
$f \cdot g = f \cdot h$ implies $g = h$ whenever the compositions make sense;
\emph{epi} or \emph{right cancellable} if
$g \cdot f = h \cdot f$ implies $g = h$ whenever the compositions make sense; and
\emph{invertible} if there is a morphism $g$ with the appropriate domain and codomain
such that $g \cdot f = \id$ and $f \cdot g = \id$.
By $\iso_\CC(A, B)$ we denote the set of all invertible
morphisms $A \to B$, and we write $A \cong B$ if $\iso_\CC(A, B) \ne \0$. Let $\Aut(A) = \iso(A, A)$.
An object $A \in \Ob(\CC)$ is \emph{rigid} if $\Aut(A) = \{\id_A\}$.

A category $\DD$ is a \emph{subcategory} of a category $\CC$ if $\Ob(\DD) \subseteq \Ob(\CC)$ and
$\hom_\DD(A, B) \subseteq \hom_\CC(A, B)$ for all $A, B \in \Ob(\DD)$.
A category $\DD$ is a \emph{full subcategory} of a category $\CC$ if $\Ob(\DD) \subseteq \Ob(\CC)$ and
$\hom_\DD(A, B) = \hom_\CC(A, B)$ for all $A, B \in \Ob(\DD)$.
If $\CC$ is a category of structures, where by a structure we mean a set together with some additional information, by
$\CC^{\fin}$ we denote the full subcategory of $\CC$ spanned by its finite members.

Let $\DD$ be a full subcategory of $\CC$.
An $S \in \Ob(\CC)$ is \emph{universal for $\DD$} if for every $D \in \Ob(\DD)$
the set $\hom_\CC(D, S)$ is nonempty and consists of monos only.
Note that if there exists an $S \in \Ob(\CC)$ universal for $\DD$ then all the morphisms
in $\DD$ are mono.
We say that $S \in \Ob(\CC)$ is \emph{projectively universal for $\DD$} if
$S$ is universal for $\DD$ in $\CC^\op$.

If $\DD$ is a full subcategory of $\CC$ the we say that $\CC$ is an \emph{ambient category} for $\DD$.
An ambient category $\CC$ is usually a category in which we can perform certain operations
that are not possible in $\DD$, or which contains an object universal for~$\DD$.

A \emph{functor} $F : \CC \to \DD$ from a category $\CC$ to a category $\DD$ maps $\Ob(\CC)$ to
$\Ob(\DD)$ and maps morphisms of $\CC$ to morphisms of $\DD$ so that
$F(f) \in \hom_\DD(F(A), F(B))$ whenever $f \in \hom_\CC(A, B)$, $F(f \cdot g) = F(f) \cdot F(g)$ whenever
$f \cdot g$ is defined, and $F(\id_A) = \id_{F(A)}$.
A functor $F : \CC \to \DD$ is an \emph{isomorphism} if there exists a functor
$E : \DD \to \CC$ such that $EF = \ID_\CC$ and $FE = \ID_\DD$, where $\ID_\CC$ denotes the identity functor
on $\CC$ which takes each object to itself and each morphism to itself. Categories $\CC$ and $\DD$ are \emph{isomorphic}
if there is an isomorphism $F : \CC \to \DD$.

A functor $U : \CC \to \DD$ is \emph{forgetful} if it is injective on morphisms in the following sense:
for all $A, B \in \Ob(\CC)$ and all $f, g \in \hom_\CC(A, B)$, if $f \ne g$ then $U(f) \ne U(g)$.
In this setting we may actually assume that $\hom_{\CC}(A, B) \subseteq \hom_\DD(U(A), U(B))$ for all $A, B \in \Ob(\CC)$.
The intuition behind this point of view is that $\CC$ is a category of structures, $\DD$ is the category of sets
and $U$ takes a structure $\calA$ to its underlying set $A$ (thus ``forgetting'' the structure). Then for every
morphism $f : \calA \to \calB$ in $\CC$ the same map is a morphism $f : A \to B$ in $\DD$.
Therefore, if $U$ is a forgetful functor we shall always take that $U(f) = f$. In particular,
$U(\id_{A}) = \id_{U(A)}$ and we, therefore, identify $\id_{A}$ with $\id_{U(A)}$.
Also, if $U : \CC \to \DD$ is a forgetful functor and all the morphisms in $\DD$ are mono, then
all the morphisms in $\CC$ are mono.

Let $F, E : \CC \to \DD$ be a pair of functors. A \emph{natural transformation from $F$ to $E$}, in symbols $\zeta : F \to E$,
is a class of arrows $\zeta_C \in \hom_\DD(F(C), E(C))$ indexed by $C \in \Ob(\CC)$ such that
\begin{center}
  \begin{tikzcd}
    F(B) \ar[r, "F(f)"] \ar[d, "\zeta_B"'] & F(C) \ar[d, "\zeta_C"] \\
    E(B) \ar[r, "E(f)"] & E(C)
  \end{tikzcd}
\end{center}
for every $B, C \in \Ob(\CC)$ and every morphism $f \in \hom_\CC(B, C)$.

\begin{EX}
  \begin{enumerate}[(1)]
  \item
    Let $\Set$ denote the category of sets and set functions and $\Set^+$ the full subcategory of $\Set$ spanned
    by all the nonempty sets.

  \item
    Let $\Top$ denote the category of topological spaces and continuous maps.

  \item
    Let $\ChEmb$ denote the category whose objects are chains and whose morphisms are embeddings.
    Let $\ChRs$ denote the category whose objects are \emph{well-ordered chains} and whose morphisms are rigid surjections.
    Let $\ChEmbFin$, resp.\ $\ChRsFin$, denote the full subcategory of $\ChEmb$, resp.\ $\ChRs$,
    spanned by finite chains.

  \item
    Let $\Alg(\Omega)$ denote the category whose objects are
    $\Omega$-algebras and morphisms are homomorphisms. Let $\AlgEpi(\Omega)$ denote the category whose objects are
    $\Omega$-algebras and morphisms are epimorphisms.

  \item
    Let $\VV$ be a variety of algebras of a fixed algebraic language. Then $\VV$ can be thought of as a category whose objects are
    the algebras in the variety and morphisms are homomorphisms. Let $\VEpi$ denote the subcategory of $\VV$
    whose objects are again all the algebras in the variety, but morphisms are epimorphisms. Finally,
    let $\VEpiFin$ be the full subcategory of $\VEpi$ spanned by its finite members.
  \end{enumerate}
\end{EX}

\paragraph{Adjunctions and monads.}
An \emph{adjunction} between categories $\BB$ and $\CC$ consists of a pair of functors $F : \BB \rightleftarrows \CC : H$
together with a family of bijections
$$
  \Phi_{X,Y} : \hom_\CC(F(X), Y) \overset\cong\longrightarrow \hom_\BB(X, H(Y))
$$
indexed by pairs $(X, Y) \in \Ob(\BB) \times \Ob(\CC)$ which is \emph{natural in both arguments} in the sense that
\begin{center}
  \begin{tikzcd}
    \hom_\CC(F(X), Y)  \arrow[r, "\Phi_{X,Y}"] \arrow[d, "(-) \cdot F(h)"'] & \hom_\BB(X,  H(Y)) \arrow[d, "(-) \cdot h"]\\
    \hom_\CC(F(X'), Y) \arrow[r, "\Phi_{X',Y}"] & \hom_\BB(X', H(Y))
  \end{tikzcd}
\end{center}
for every $h \in \hom_\BB(X', X)$, and
\begin{center}
  \begin{tikzcd}
    \hom_\CC(F(X), Y)  \arrow[d, "g \cdot (-)"'] \arrow[r, "\Phi_{X,Y}"] & \hom_\BB(X, H(Y)) \arrow[d, "F(g) \cdot (-)"] \\
    \hom_\CC(F(X), Y') \arrow[r, "\Phi_{X,Y'}"] & \hom_\BB(X, H(Y'))
  \end{tikzcd}
\end{center}
every $g \in \hom_\CC(Y, Y')$. The functor $F$ is then \emph{left adjoint} (to $H$) and $H$ is \emph{right adjoint} (to $F$).

Let $\CC$ be a category and $T : \CC \to \CC$ a functor. \emph{Multiplication for $T$} is a natural transformation
$\mu : TT \to T$ such that for each $A \in \Ob(\CC)$:
\begin{center}
    \begin{tikzcd}
      TTT(A) \arrow[r, "T(\mu_A)"] \arrow[d, "\mu_{T(A)}"'] & TT(A) \arrow[d, "\mu_A"] \\
      TT(A) \arrow[r, "\mu_A"'] & T(A)
    \end{tikzcd}
\end{center}
A natural transformation $\eta : \ID \to T$ is a \emph{unit for $\mu$} if
\begin{center}
    \begin{tikzcd}
      T(A) \arrow[r, "T(\eta_A)"] \arrow[dr, "\id_{T(A)}"'] & TT(A) \arrow[d, "\mu_A"] & T(A) \arrow[l, "\eta_{T(A)}"'] \arrow[dl, "\id_{T(A)}"] \\
         & T(A) &
    \end{tikzcd}
\end{center}
\noindent
for each $A \in \Ob(\CC)$. A \emph{monad} on a category $\CC$ is a triple $(T, \mu, \eta)$ where $T : \CC \to \CC$ is a functor,
$\mu$ is a multiplication for $T$ and $\eta$ is a unit for~$\mu$.

Let $F : \CC \to \CC$ be a functor. An \emph{$F$-algebra} is a pair $(A, \alpha)$ where $\alpha \in \hom_\CC(F(A), A)$, while
an \emph{$F$-coalgebra} is a pair $(A, \alpha)$ where $\alpha \in \hom_\CC(A, F(A))$.
An \emph{algebraic homomorphism} between $F$-algebras $(A, \alpha)$ and $(B, \beta)$ is a morphism $f \in \hom_\CC(A, B)$ such that
the diagram below commutes:
\begin{center}
    \begin{tikzcd}
      F(A) \arrow[r, "F(f)"] \arrow[d, "\alpha"'] & F(B) \arrow[d, "\beta"]\\
      A \arrow[r, "f"'] & B
    \end{tikzcd}
\end{center}

Two categories are traditionally associated to each monad $(T, \mu, \eta)$: the Kleisli category $\KK = \KK(T, \mu, \eta)$
and the Eilenberg-Moore category $\EM = \EM(T, \mu, \eta)$.
The objects of the Kleisli category $\KK(T, \mu, \eta)$ are the same as the objects of $\CC$, morphisms are defined by
$$
  \hom_\KK(A, B) = \hom_\CC(A, T(B))
$$
and the composition in $\KK$ for $f \in \hom_\KK(A, B)$ and $g \in \hom_\KK(B, C)$ is given by
$$
  g \mathbin{\cdot_\KK} f = \mu_C \cdot T(g) \cdot f.
$$
The objects of the Eilenberg-Moore category $\EM(T, \mu, \eta)$ are Eilenberg-Moore $T$-algebras (special $T$-algebras
to be defined immediately), morphisms are algebraic homomorphisms and the composition is as in $\CC$.
An \emph{Eilenberg-Moore $T$-algebra} is a $T$-algebra for which the following two diagrams commute:
\begin{equation}\label{rpemklei.eq.def-Talg}
  \begin{array}{c@{\qquad}c}
    \begin{tikzcd}
      TT(A) \arrow[r, "T(\alpha)"] \arrow[d, "\mu_A"'] & T(A) \arrow[d, "\alpha"] \\
      T(A) \arrow[r, "\alpha"'] & A
    \end{tikzcd}
    &
    \begin{tikzcd}
      A \arrow[r, "\eta_A"] \arrow[dr, "\id_A"'] & T(A) \arrow[d, "\alpha"] \\
      & A
    \end{tikzcd}
  \end{array}
\end{equation}
A \emph{weak Eilenberg-Moore $T$-algebra} is a $T$-algebra for which only the diagram on the left commutes.
Let $\EM^w(T, \mu)$ denote the category of weak Eilenberg-Moore $T$-algebras and algebraic homomorphisms.
A \emph{free Eilenberg-Moore $T$-algebra} is an Eilenberg-Moore $T$-algebra of the form $(T(A), \mu_A)$.

\begin{EX}\label{rpemklei.ex.T-Sigma}
  Let $\Omega$ be an algebraic language, that is, a set of functional and constant symbols.
  For any nonempty set $X$ let $T(X)$ denote the set of all $\Omega$-terms in variables from~$X$.
  For a function $f : X \to Y$ define $T(f) : T(X) \to T(Y)$ to be the substitution of variables
  with respect to~$f$. Then $T : \Set^+ \to \Set^+$ is a functor.
  Let $\eta_X : X \to T(X)$ send $x$ as a variable to $x$ as an $\Omega$-term
  and let $\mu_X : T T(X) \to T(X)$ denote the substitution of terms for variables. Then
  $(T, \mu, \eta)$ is a $\Set^+$-monad such that the category $\EM(T, \mu, \eta)$ is isomorphic to $\Alg(\Omega)$.
  Let us outline briefly the way the isomorphism between the two categories is constructed.

  For an $\Omega$-algebra $\calA = (A, \Omega^A)$ the corresponding $T$-algebra $(A, \eval^A)$ is
  constructed so that the structure map $\eval^A : T(A) \to A$ is just the evaluation in~$\calA$.

  Conversely, every $T$-algebra $(A, \eval^A)$ corresponds to a unique $\Omega$-algebra $\calA = (A, \Omega^A)$.
  To see that this is indeed the case, let us show how to decode the fundamental operations of $\calA$ from $(A, \eval^A)$.
  Take a function symbol $f \in \Omega$, say of arity $n$, and define the operation $f^A : A^n \to A$ by
  $$
    f^A(a_1, \ldots, a_n) = \eval^A(f(\eta(a_1), \ldots, \eta(a_n))),
  $$
  where $f(\eta(a_1), \ldots, \eta(a_n)) \in T(A)$ is a term.
  Note that the same procedure applies to constant symbols as well, since they can be (and usually are)
  treated as 0-ary function symbols. We thus get an $\Omega$-algebra $\calA = (A, \Omega^A)$.
  The requirements \eqref{rpemklei.eq.def-Talg} then ensure that $(A, \eval^A)$
  is the $T$-algebra that corresponds to $\calA$.

  The two categories are isomorphic simply because we are given two different presentations
  of the same object: an $\Omega$-algebra $\calA = (A, \Omega^A)$ can be thought of as a ``small'' presentation of the
  way things are computed in $\calA$ as it provides the ``multiplication tables'' of fundamental operations only.
  On the other hand, the corresponding $T$-algebra $(A, \eval^A)$ is defined on the same set, but the evaluation map
  $\eval^A : T(A) \to A$ provides the value of each possible expression that can be formed from the elements of $A$ and
  fundamental operations in $\Omega$. See~\cite[Sections 10.3 and 10.5]{awodey} for further details and concrete examples.
\end{EX}

\section{Ramsey properties in a category}
\label{rpemklei.sec.RP}

In this section we collect and prove several results about the Ramsey property, Ramsey degrees,
dual Ramsey property and dual small Ramsey degrees in a category. We then use the results of this section
as the main tool to obtain new Ramsey results about algebras in a variety (in Section~\ref{rpemklei.sec.drp-V}).

\paragraph{The arrow notation.}
For $k \in \NN$, a $k$-coloring of a set $S$ is any mapping $\chi : S \to k$, where,
as usual, we identify $k$ with $\{0, 1,\ldots, k-1\}$.
Let $\CC$ be a locally small category.
For integers $k \ge 2$ and $t \ge 1$, and objects $A, B, C \in \Ob(\CC)$ we write
$$
  C \longrightarrow (B)^{A}_{k, t}
$$
to denote that for every $k$-coloring $\chi : \hom(A, C) \to k$ there is a morphism
$w \in \hom(B, C)$ such that $|\chi(w \cdot \hom(A, B))| \le t$.
(For a set of morphisms $F$ we let $w \cdot F = \{ w \cdot f : f \in F \}$.)
In case $t = 1$ we write
$
  C \longrightarrow (B)^{A}_{k}.
$
We write
$$
  C \longleftarrow (B)^{A}_{k, t}, \text{ resp.\ } C \longleftarrow (B)^{A}_{k},
$$
to denote that $C \longrightarrow (B)^{A}_{k, t}$, resp.\ $C \longrightarrow (B)^{A}_{k}$, in $\CC^\op$.

\begin{LEM}\label{rpemklei.lem.C-D} \cite[Lemma 2.4]{masul-scow}
  Let $\CC$ be a locally small category such that all the morphisms in $\CC$ are mono and let $A, B, C, D \in \Ob(\CC)$. If
  $C \longrightarrow (B)^A_{k,t}$ for some $k, t \ge 2$ and if $C \overset\CC\longrightarrow D$, then $D \longrightarrow (B)^A_{k,t}$.\qed
\end{LEM}

The above lemma tells us that with the arrow notation we can always go ``up to a superstructure'' of $C$.
In some cases we can also ``go down to a substructure'' of $C$.
Let $\AAA$ be a subcategory of $\CC$ and let $C \in \Ob(\CC)$. An object $B \in \Ob(\AAA)$ together with a
morphism $c : B \to C$ is a \emph{coreflection of $C$ in $\AAA$} if for every $A \in \Ob(\AAA)$ and every morphism
$f \in \hom_\CC(A, C)$ there is a unique morphism $g \in \hom_\AAA(A, B)$ such that $c \cdot g = f$:
\begin{center}
  \begin{tikzcd}
         &    C & & \CC \\
    A \arrow[r, "g"'] \arrow[ur,"f"] & B \arrow[u, "c"'] & & \AAA 
  \end{tikzcd}
\end{center}

Dually, an object $B \in \Ob(\AAA)$ together with a
morphism $r : C \to B$ is a \emph{reflection of $C$ in $\AAA$} if for every $A \in \Ob(\AAA)$ and every morphism
$f \in \hom_\CC(C, A)$ there is a unique morphism $g \in \hom_\AAA(B, A)$ such that $g \cdot r = f$:
\begin{center}
  \begin{tikzcd}
         &    C \arrow[d, "r"] \arrow[dl, "f"'] & & \CC \\
    A & B \arrow[l, "g"] & & \AAA
  \end{tikzcd}
\end{center}

\begin{LEM}\label{rpemklei.lem.refl}
  Let $\CC$ be a locally small category such that all the morphisms in $\CC$ are mono. Let $\AAA$ be a full
  subcategory of $\CC$, let $A, B, D \in \Ob(\AAA)$ and $C \in \Ob(\CC)$. If
  $C \longrightarrow (B)^A_{k,t}$ for some $k, t \ge 2$ and if
  $c : D \to C$ is a coreflection of $C$ in $\AAA$ then $D \longrightarrow (B)^A_{k,t}$.

  Dually, let $\CC$ be a locally small category such that all the morphisms in $\CC$ are epi. Let $\AAA$ be a full
  subcategory of $\CC$, let $A, B, D \in \Ob(\AAA)$ and $C \in \Ob(\CC)$. If
  $C \longleftarrow (B)^A_{k,t}$ for some $k, t \ge 2$ and if
  $r : C \to D$ is a reflection of $C$ in $\AAA$ then $D \longleftarrow (B)^A_{k,t}$.
\end{LEM}
\begin{proof}
  Take any coloring $\chi : \hom(A, D) \to k$ and define $\chi' : \hom(A, C) \to k$ as follows:
  $\chi'(c \cdot g) = \chi(g)$ for all $g \in \hom(A, D)$, and $\chi'(f) = 0$ for all $f \in \hom(A, C) \setminus c \cdot \hom(A, D)$.
  Note that the definition of $\chi'$ is correct because $c$ is mono. Since $C \longrightarrow (B)^A_{k,t}$
  there is a $w' \in \hom(B, C)$ such that
  $$
    |\chi'(w' \cdot \hom(A, B))| \le t.
  $$
  Because $c : D \to C$ is a coreflection of $C$ in $\AAA$ there is a unique morphism $w \in \hom(B, D)$ such that $c \cdot w = w'$:
\begin{center}
  \begin{tikzcd}
         &    C                           &   \\
    A \arrow[r, "g"'] & D \arrow[u, "c"] & B \arrow[ul, "w'"'] \arrow[l, "w"] 
  \end{tikzcd}
\end{center}
  So,
  $
    |\chi'(c \cdot w \cdot \hom(A, B))| \le t
  $
  whence
  $
    |\chi(w \cdot \hom(A, B))| \le t
  $.
\end{proof}

\paragraph{Finite Ramsey phenomena.}
A category $\CC$ has the \emph{(finite) Ramsey property} if
for every integer $k \ge 2$ and all $A, B \in \Ob(\CC)$ there is a
$C \in \Ob(\CC)$ such that $C \longrightarrow (B)^{A}_k$.
A category $\CC$ has the \emph{(finite) dual Ramsey property} if $\CC^\op$ has the Ramsey property.

Both the Ramsey property and the dual Ramsey property impose severe restrictions on the classes of objects enjoying the property.
For example, all the objects in such classes have to be rigid:

\begin{THM} (\cite[Proposition 2.3]{masul-scow}, cf.~\cite{Nesetril,zucker1})\label{rpemklei---rigid}
  Let $\CC$ be a locally small category such that all the morphisms in $\CC$ are mono. If $\CC$ has the Ramsey property then all the
  objects in $\CC$ are rigid.
  
  Dually, let $\CC$ be a locally small category such that all the morphisms in $\CC$ are epi. If $\CC$ has the dual Ramsey property then all the
  objects in $\CC$ are rigid.\qed
\end{THM}

For $A \in \Ob(\CC)$ let $t_\CC(A)$ denote the least positive integer $n$ such that
for all $k \ge 2$ and all $B \in \Ob(\CC)$ there exists a $C \in \Ob(\CC)$ such that
$C \longrightarrow (B)^{A}_{k, n}$, if such an integer exists.
Otherwise put $t_\CC(A) = \infty$. The number $t_\CC(A)$ is referred to as the \emph{small Ramsey degree}
of $A$ in $\CC$. A category $\CC$ has the \emph{finite small Ramsey degrees} if $t_\CC(A) < \infty$
for all $A \in \Ob(\CC)$. Clearly, a category $\CC$ has the Ramsey property if and only if $t_\CC(A) = 1$
for all $A \in \Ob(\CC)$. In this parlance the Finite Ramsey Theorem takes the following form.

\begin{THM} [The Finite Ramsey Theorem~\cite{ramsey}]
  The category $\ChEmbFin$ has the Ramsey property.\qed
\end{THM}

By straightforward dualization we can introduce dual small Ramsey degrees
$t^\partial_\CC(A)$ by $t^\partial_\CC(A) = t_{\CC^\op}(A)$.
We then say that a category $\CC$ has the \emph{finite dual small Ramsey degrees} if $\CC^\op$ has finite small Ramsey degrees.
Clearly, a category $\CC$ has the dual Ramsey property if and only if $t^\partial_\CC(A) = 1$
for all $A \in \Ob(\CC)$. The Finite Dual Ramsey Theorem~\cite{graham-rothschild-DRT} takes the following form.

\begin{THM}[The Finite Dual Ramsey Theorem~\cite{graham-rothschild-DRT}]
  The category $\ChRsFin$ has the dual Ramsey property.\qed
\end{THM}

The Ramsey property for ordered structures implies the existence of finite small Ramsey degrees for
the corresponding unordered structures. This was first observed for categories of structures
in~\cite{dasilvabarbosa}, and generalized to arbitrary categories in~\cite{masul-dual-kpt}.
This generalization will prove useful here because it will enable us to derive statements
about the dual small Ramsey degrees for finite algebras in a variety.

Let us outline the main tool we employ to obtain results of this form.
Following \cite{KPT,dasilvabarbosa,masul-dual-kpt} we say that
an \emph{expansion} of a category $\CC$ is a category $\CC^*$ together with a forgetful functor $U : \CC^* \to \CC$.
We shall generally follow the convention that $A, B, C, \ldots$ denote objects from $\CC$
while $\calA, \calB, \calC, \ldots$ denote objects from $\CC^*$.
Since $U$ is injective on hom-sets we may safely assume that
$\hom_{\CC^*}(\calA, \calB) \subseteq \hom_\CC(A, B)$ where $A = U(\calA)$, $B = U(\calB)$.
In particular, $\id_\calA = \id_A$ for $A = U(\calA)$. Moreover, it is safe to drop
subscripts $\CC$ and $\CC^*$ in $\hom_\CC(A, B)$ and $\hom_{\CC^*}(\calA, \calB)$, so we shall
simply write $\hom(A, B)$ and $\hom(\calA, \calB)$, respectively.
Let
$
  U^{-1}(A) = \{\calA \in \Ob(\CC^*) : U(\calA) = A \}
$. Note that this is not necessarily a set.

An expansion $U : \CC^* \to \CC$ is \emph{reasonable} (cf.~\cite{KPT,masul-dual-kpt}) if
for every $e \in \hom(A, B)$ and every $\calA \in U^{-1}(A)$ there is a $\calB \in U^{-1}(B)$ such that
$e \in \hom(\calA, \calB)$:
\begin{center}
    \begin{tikzcd}
      \calA \arrow[r, "e"] \arrow[d, dashed, mapsto, "U"'] & \calB \arrow[d, dashed, mapsto, "U"] \\
      A \arrow[r, "e"'] & B
    \end{tikzcd}
\end{center}
An expansion $U : \CC^* \to \CC$ has \emph{unique restrictions} \cite{masul-dual-kpt} if
for every $\calB \in \Ob(\CC^*)$ and every $e \in \hom(A, U(\calB))$ there is a \emph{unique} $\calA \in U^{-1}(A)$
such that $e \in \hom(\calA, \calB)$:
\begin{center}
    \begin{tikzcd}
      \llap{\hbox{$\restr \calB e = \mathstrut$}}\calA \arrow[r, "e"] \arrow[d, dashed, mapsto, "U"'] & \calB \arrow[d, dashed, mapsto, "U"] \\
      A \arrow[r, "e"'] & B
    \end{tikzcd}
\end{center}
We denote this unique $\calA$ by $\restr \calB e$ and refer to it as the \emph{restriction of $\calB$ along~$e$}.

The following result was first proved for categories of structures in \cite{dasilvabarbosa},
and for general categories in \cite{masul-dual-kpt}. This more general statement will
be of use in Section~\ref{rpemklei.sec.drp-V} where we shall need the dual of the statement.

\begin{THM}\label{sbrd.thm.small1} \cite{dasilvabarbosa,masul-dual-kpt}
  Let $\CC$ and $\CC^*$ be locally small categories such that all the morphisms in $\CC$ and $\CC^*$ are mono.
  Let $U : \CC^* \to \CC$ be a reasonable expansion with unique restrictions such that $U^{-1}(A)$ is a set for all $A \in \Ob(\CC)$.
  For any $A \in \Ob(\CC)$ we then have:
  $$
    t_{\CC}(A) \le \sum_{\calA \in U^{-1}(A)} t_{\CC^*}(\calA).
  $$
  Consequently, if $U^{-1}(A)$ is finite and
  $t_{\CC^*}(\calA) < \infty$ for all $\calA \in U^{-1}(A)$ then $t_\CC(A) < \infty$.

  In particular, if $U : \CC^* \to \CC$ is a reasonable expansion with unique restrictions such that $\CC^*$ has the
  Ramsey property and $U^{-1}(A)$ is finite for all $A \in \Ob(\CC)$ then $\CC$ has finite small Ramsey degrees.\qed
\end{THM}

\paragraph{Infinite Ramsey phenomena.}
Let $\CC$ be a locally small category.
For $A, S \in \Ob(\CC)$ let $T_\CC(A, S)$ denote the least positive integer $n$ such that
$S \longrightarrow (S)^{A}_{k, n}$ for all $k \ge 2$, if such an integer exists.
Otherwise put $T_\CC(A, S) = \infty$. The number $T_\CC(A, S)$ is referred to as the \emph{big Ramsey degree}
of $A$ in $S$.
By straightforward dualization, we can introduce \emph{dual big Ramsey degrees}
$T^\partial_\CC(A, S)$ by $T^\partial_\CC(A, S) = T_{\CC^\op}(A, S)$.
We shall drop the category in the index whenever it is clearly stated which category we work in.
In this parlance the Infinite Ramsey Theorem takes the following form.

\begin{THM}[The Infinite Ramsey Theorem~\cite{ramsey}] \label{rpemklei.thm.IRT}
  In the category $\ChEmb$ we have that
  $T(A, \omega) = 1$ for every finite chain~$A$.\qed
\end{THM}

The Infinite Dual Ramsey Theorem of Carlson and Simpson \cite{carlson-simpson-1984}
requires additional infrastructure and the notion of
dual big Ramsey degrees with respect to special colorings.

A category $\CC$ is \emph{enriched over $\Top$} if each $\hom_\CC(A, B)$ is a topological space
and the composition $\Boxed\cdot : \hom_\CC(B, C) \times \hom_\CC(A, B) \to \hom_\CC(A, C)$ is continuous
for all $A, B, C \in \Ob(\CC)$. Any locally small category can be thought of as a category enriched over $\Top$ by taking each
hom-set to be a discrete space. We shall refer to this as the \emph{discrete enrichment}.
(Note that a category enriched over $\Top$ has to be locally small.)

  The category $\ChRs$ can be enriched over $\Top$ in a nontrivial way as follows:
  each hom-set $\hom_{\ChRs}((A, \Boxed<), (B, \Boxed<))$ inherits the topology from the Tychonoff topology on $B^A$
  with both $A$ and $B$ discrete. Whenever we refer to $\ChRs$ as a category enriched over $\Top$ we have this particular
  enrichment in mind.

For a topological space $X$ and an integer $k \ge 2$ a \emph{Borel $k$-coloring of $X$} is any mapping $\chi : X \to k$
such that $\chi^{-1}(i)$ is a Borel set for all~$i \in k$.

Let $\CC$ be a category enriched over $\Top$.
For $A, B, C \in \Ob(\CC)$ we write $C \Borel\longrightarrow (B)^{A}_{k, n}$
to denote that for every Borel $k$-coloring $\chi : \hom_\CC(A, C) \to k$
there is a morphism $w \in \hom_\CC(B, C)$ such that $|\chi(w \cdot \hom_\CC(A, B))| \le n$.
For $A, S \in \Ob(\CC)$ let $\flatT_\CC(A, S)$ denote the least positive integer $n$ such that
$S \Borel\longrightarrow (S)^{A}_{k, n}$ for all $k \ge 2$, if such an integer exists.
Otherwise put $\flatT_\CC(A, S) = \infty$. The number $\flatT_\CC(A, S)$ will be referred to as the
\emph{big Ramsey degree of $A$ in $S$ with respect to Borel colorings}.
By straightforward dualization, we can introduce \emph{big dual Ramsey degrees with respect to Borel colorings}
$\flatdualT_\CC(A, S)$ by $\flatdualT_\CC(A, S) = \flatT_{\CC^\op}(A, S)$.
The Infinite Dual Ramsey Theorem of Carlson and Simpson~\cite{carlson-simpson-1984}
now takes the following form.

\begin{THM}[The Infinite Dual Ramsey Theorem~\cite{carlson-simpson-1984}]\label{rpemklei.thm.IDRT}
  In the category $\ChRs$ enriched over $\Top$ as above,
  $\flatdualT(A, \omega) = 1$ for every finite chain~$A$.\qed
\end{THM}

The following result was first proved for categories of structures in \cite{dasilvabarbosa},
and for general categories in \cite{masul-big-v-small}.

\begin{THM}\label{rpemklei---big.v.small} (cf. \cite{dasilvabarbosa,masul-big-v-small})
  Let $\CC$ and $\CC^*$ be locally small categories such that all the morphisms in $\CC$ and $\CC^*$ are mono.
  Let $U : \CC^* \to \CC$ be an expansion with unique restrictions.
  For $A \in \Ob(\CC)$, $S^* \in \Ob(\CC^*)$ and $S = U(S^*)$, if $U^{-1}(A)$ is finite then
  $$
    T_{\CC}(A, S) \le \sum_{A^* \in U^{-1}(A)} T_{\CC^*}(A^*, S^*).\qed
  $$
\end{THM}

However, for the results in Section~\ref{rpemklei.sec.drp-V} we shall need the
version of this statement which, with Theorem~\ref{rpemklei.thm.IDRT} in mind, has to be formulated
in terms of categories enriched over $\Top$.
Let $\BB$ and $\CC$ be categories enriched over $\Top$.

\begin{THM}\label{sbrd.thm.big1-dual-borel}
  Let $\CC$ and $\CC^*$ be categories enriched over $\Top$ such that morphisms in both $\CC$ and $\CC^*$ are mono.
  Let $U : \CC^* \to \CC$ be an expansion with unique restrictions and assume that $\hom_{\CC^*}(A^*, B^*)$ is a Borel
  set in $\hom_{\CC}(U(A^*), U(B^*))$ for all $A^*, B^* \in \Ob(\CC^*)$.
  For $A \in \Ob(\CC)$, $S^* \in \Ob(\CC^*)$ and $S = U(S^*)$, if $U^{-1}(A)$ is finite then
  $$
    \flatT_{\CC}(A, S) \le \sum_{A^* \in U^{-1}(A)} \flatT_{\CC^*}(A^*, S^*).
  $$
\end{THM}
\begin{proof}
  The proof is analogous to the proof of the corresponding
  statement in~\cite{masul-big-v-small}. The only difference is that we now have to ensure that
  all the colorings we construct along the way are Borel. But that follows immediately from the 
  additional assumption that $\hom_{\CC^*}(A^*, B^*)$ is a Borel
  set in $\hom_{\CC}(U(A^*), U(B^*))$ for all $A^*, B^* \in \Ob(\CC^*)$,
  and the fact that in categories enriched over $\Top$ the composition of morphisms is continuous.
\end{proof}

\paragraph{Ramsey properties and adjunctions.}
Our major tool for transporting the Ramsey property from one context to another
is to establish an adjunction-like relationship between the corresponding categories.

\begin{THM}\label{rpemklei.thm.adjoints-rp} \cite{masul-scow}
  Right adjoints preserve the Ramsey property while left adjoints preserve the dual Ramsey property.
  More precisely, let $\BB$ and $\CC$ be locally small categories and let $F : \BB \rightleftarrows \CC : H$ be an adjunction.
  
  $(a)$ If $\CC$ has the Ramsey property then so does $\BB$.
  
  $(b)$ If $\BB$ has the dual Ramsey property then so does $\CC$.\qed
\end{THM}

\begin{THM}
  Let $\CC$ be a locally small category, $(T, \mu, \eta)$ a monad on $\CC$, and let
  $\KK = \KK(T, \mu, \eta)$ and $\EM = \EM(T, \mu, \eta)$
  be the Kleisli category and the Eilenberg-Moore category, respectively, for the monad. If $\CC$ has the
  dual Ramsey property then so do both $\KK$ and $\EM$.
\end{THM}
\begin{proof}
  Let $(T, \mu, \eta)$ be a monad on $\CC$, and let $\KK = \KK(T, \mu, \eta)$ and $\EM = \EM(T, \mu, \eta)$
  be the Kleisli category and the Eilenberg-Moore category, respectively, for the monad.
  It is a well-known fact (see~\cite{maclane}) that there exist adjunctions
  $\CC \rightleftarrows \KK$ and $\CC \rightleftarrows \EM$.
  The statement now follows from Theorem~\ref{rpemklei.thm.adjoints-rp}.
\end{proof}

However, more is true in case of the Eilenberg-Moore construction.
We are now going to show that the (dual) Ramsey property
carries over from $\CC$ to a more general context of (co)algebras for functors with (co)multiplication,
which are straightforward weakenings of (co)monads. The proof relies on the following weakening of the notion of adjunction.

\begin{DEF}\label{opos.def.PA} \cite{masul-preadj}
  Let $\BB$ and $\CC$ be locally small categories. A pair of maps
  $
    F : \Ob(\BB) \rightleftarrows \Ob(\CC) : H
  $
  is a \emph{pre-adjunction between $\BB$ and $\CC$} provided there is a family of maps
  $
    \Phi_{X,Y} : \hom_\CC(F(X), Y) \to \hom_\BB(X, H(Y))
  $
  indexed by the pairs $(X, Y) \in \Ob(\BB) \times \Ob(\CC)$ and satisfying the following:
  \begin{itemize}
  \item[(PA)]
  for every $C \in \Ob(\CC)$, every $A, B \in \Ob(\BB)$,
  every $u \in \hom_\CC(F(B), C)$ and every $f \in \hom_\BB(A, B)$ there is a $v \in \hom_\CC(F(A), F(B))$
  satisfying $\Phi_{B, C}(u) \cdot f = \Phi_{A, C}(u \cdot v)$.
  \end{itemize}
  \begin{center}
    \begin{tikzcd}
        F(B) \arrow[rr, "u"]                          &     & C & & B \arrow[rr, "\Phi_{B, C}(u)"]                          &      & H(C) \\
        F(A) \arrow[u, "v"] \arrow[urr, "u \cdot v"'] &     &       & & A \arrow[u, "f"] \arrow[urr, "\Phi_{A, C}(u \cdot v)"'] \\
                                                            & \CC \arrow[rrrr, shift left, "H"]
                                                            &       & &                                                                      & \BB \arrow[llll, shift left, "F"]
    \end{tikzcd}
  \end{center}
\end{DEF}
\noindent
(Note that in a pre-adjunction $F$ and $H$ are \emph{not} required to be functors, just maps from the class of objects of one of the two
categories into the class of objects of the other category; also $\Phi$ is not required to be a natural isomorphism, just a family of
maps between hom-sets satisfying the requirement above.)

\begin{THM}\label{opos.thm.main} \cite{masul-preadj}
  Let $\BB$ and $\CC$ be locally small categories and let $F : \Ob(\BB) \rightleftarrows \Ob(\CC) : H$ be a pre-adjunction.
  \begin{itemize}
  \item[$(a)$]
    Let $t, k \ge 2$ be integers, let $A, B \in \Ob(\BB)$ and $C \in \Ob(\CC)$.
    If $C \longrightarrow (F(B))^{F(A)}_{k,t}$ in $\CC$ then $H(C) \longrightarrow (B)^A_{k,t}$ in $\BB$.
  \item[$(b)$]
    $t_\BB(A) \le t_\CC(F(A))$ for all $A \in \Ob(\BB)$.
  \item[$(c)$]
    If $\CC$ has the Ramsey property then so does $\BB$.\qed
  \end{itemize}  
\end{THM}

Although pre-adjunctions constitute relatively loose relationships between categories,
they turn out to be a useful strategy of transporting the Ramsey property from one category to the other.
This strategy is motivated by the proof of \cite[Theorem 12.13]{promel-book}
where the Ramsey property for finite ordered graphs was shown to be a direct
consequence of the Graham–Rothschild Theorem. It was then used in~\cite{masul-preadj}
to provide new proofs of the Ramsey property for the category
of all finite posets with linear extensions (by establishing a pre-adjunction with the Graham–Rothschild
category), for the category of finite convexly ordered ultrametric spaces (by establishing a pre-adjunction with the
category of finite posets with linear extensions) and the category of finite linearly ordered metric spaces
(by establishing a different pre-adjunction with the category of finite posets with linear extensions).

\begin{THM}\label{rpemklei.thm.ramsey}
  Let $\CC$ be a locally small category, $T : \CC \to \CC$ a functor and $\mu : TT \to T$ a multiplication for $T$.
  Let $\WW$ be a category whose objects are weak Eilenberg-Moore $T$-algebras, morphisms are algebraic
  homomorphisms and the composition of morphisms is as in $\CC$. If $\CC$ has the dual Ramsey property then so does
  every full subcategory of $\WW$ which contains all the free Eilenberg-Moore $T$-algebras.
\end{THM}
\begin{proof}
  Let $\BB$ be a full subcategory of $\WW$ such that
  all the free Eilenberg-Moore $T$-algebras are in~$\BB$.
  By Theorem~\ref{opos.thm.main} in order to show that $\BB$ has the dual Ramsey property
  it suffices to construct a pre-adjunction $F : \Ob(\BB^\op) \rightleftarrows \Ob(\CC^\op) : H$.
  For $(B, \beta) \in \Ob(\BB)$ put $F(B, \beta) = B$,
  for $C \in \Ob(\CC)$ put $H(C) = (T(C), \mu_C)$ and
  for $u \in \hom_\CC(C, B)$ put $\Phi_{(B, \beta), C}(u) = \beta \cdot T(u)$.

  Let us first show that the definition of $\Phi$ is correct by showing that for each
  $u \in \hom_\CC(C, B)$ we have that $\Phi_{(B, \beta), C}(u)$ is a homomorphism
  from $H(C)$ to $(B, \beta)$, that is:
  \begin{center}
    \begin{tikzcd}
      TT(C) \arrow[d, "\mu_C"'] \arrow[rr, "T(\beta \cdot T(u))"] & & T(B) \arrow[d, "\beta"] \\
      T(C) \arrow[rr, "\beta \cdot T(u)"'] & & B
    \end{tikzcd}
  \end{center}
  \noindent
  But this is straightforward -- the two squares below commute because $\mu : TT \to T$ is a natural transformation
  and $(B, \beta)$ is a weak $T$-algebra:
  \begin{center}
      \begin{tikzcd}
        TT(C) \arrow[r, "TT(u)"] \arrow[d, "\mu_C"'] & TT(B) \arrow[d, "\mu_B"] \arrow[r, "T(\beta)"]  & T(B) \arrow[d, "\beta"]\\
        T(C) \arrow[r, "T(u)"'] & T(B) \arrow[r, "\beta"'] & B
      \end{tikzcd}
  \end{center}
  \noindent
  To complete the proof we still have to show that the condition (PA) of Definition~\ref{opos.def.PA}
  is satisfied.
  Take any $C \in \Ob(\CC)$ and $(A, \alpha), (B, \beta) \in \Ob(\BB)$, take
  arbitrary $u \in \hom_\CC(C, B)$ and an
  arbitrary homomorphism $f \in \hom_\BB((B, \beta), (A, \alpha))$.
  Take $v = f \in \hom_\CC(B, A)$ and note that
  $$
    f \cdot \Phi_{(B, \beta), C}(u)
    = f \cdot \beta \cdot T(u)
    = \alpha \cdot T(f) \cdot T(u)
    = \alpha \cdot T(f \cdot u)  = \Phi_{(A, \alpha), C}(f \cdot u)
  $$
  having in mind that $f \cdot \beta = \alpha \cdot T(f)$ because $f$ is a homomorphism.
  This completes the proof.
\end{proof}

It might be of interest to note that the pre-adjunction constructed in the proof of Theorem~\ref{rpemklei.thm.ramsey}
cannot be expected to be an adjunction. Although $F$ and $H$ behave as the forgetful functor and the free-algebra functor,
respectively, which is usual in adjunctions associated with categories of algebras, one cannot hope for $\Phi$ to be either
bijective on hom-sets or natural in its arguments.

\section{Dual Ramsey properties for varieties of algebras}
\label{rpemklei.sec.drp-V}

This section is devoted to the study of dual Ramsey phenomena in nontrivial varieties of algebras.
Motivated by the fact that the category of weak Eilenberg-Moore algebras for a monad
has the dual Ramsey property (Theorem~\ref{rpemklei.thm.ramsey}), we show
that for every algebraic language $\Omega$ and every nontrivial variety $\VV$ of $\Omega$-algebras
the class of finite ordered $\VV$ algebras taken with rigid epimorphisms (that is, epimorphisms of algebras that are at the
same time rigid surjections) has the dual Ramsey property. The unordered version then follows
immediately: for every algebraic language $\Omega$ and every nontrivial variety $\VV$ of $\Omega$-algebras,
finite $\VV$ algebras have finite small dual Ramsey degrees with respect to epimorphisms.

Our proof mimics the proof of the fact that the category of weak Eilenberg-Moore algebras for a monad
has the dual Ramsey property (Theorem~\ref{rpemklei.thm.ramsey}). Unfortunately,
we are unable to apply Theorem~\ref{rpemklei.thm.ramsey} directly
because the monad we are going to construct will not produce finite free algebras.
We shall bypass this issue with the help of
the following compactness argument which was first proved for categories of structures in~\cite{mu-pon} (see also~\cite{vanthe-more}),
and for general categories in \cite{masul-dual-kpt}.

Let $\BB$ be a full subcategory of a locally small category $\CC$.
An $F \in \Ob(\CC)$ is \emph{weakly locally finite for $\BB$} (cf.\ locally finite in~\cite{masul-dual-kpt}) if
for every $A, B \in \Ob(\BB)$ and every $e \in \hom(A, F)$, $f \in \hom(B, F)$ there exist $D \in \Ob(\BB)$,
$r \in \hom(D, F)$, $p \in \hom(A, D)$ and $q \in \hom(B, D)$ such that $r \cdot p = e$ and $r \cdot q = f$.
We say that $F \in \Ob(\CC)$ is \emph{projectively weakly locally finite for $\BB$} if
$F$ is weakly locally finite for $\BB$ in $\CC^\op$.

\begin{LEM}\label{akpt.lem.ramseyF} \cite{mu-pon,vanthe-more,masul-dual-kpt}
  Let $\BB$ be a full subcategory of a locally small category $\CC$ and fix an $F \in \Ob(\CC)$ which is universal and weakly locally finite
  for $\BB$. The following are equivalent for all $t \ge 2$ and all $A \in \Ob(\BB)$:
  \begin{enumerate}\def\labelenumi{(\arabic{enumi})}
  \item
    $t_{\BB}(A) \le t$;
  \item
    $F \longrightarrow (B)^A_{k, t}$ for all $k \ge 2$ and all $B \in \Ob(\BB)$ such that $A \overset\BB\longrightarrow B$.\qed
  \end{enumerate}
\end{LEM}

We then prove that for every countable algebraic language $\Omega$ and every nontrivial variety $\VV$ of $\Omega$-algebras,
finite ordered $\VV$ algebras have finite big Ramsey degrees in the ordered free $\VV$ algebra
$\hat \calF_\VV(\omega)$ on $\omega$ generators with respect to Borel colorings.
The corresponding result for the unordered case follows by a straightforward modification of the ideas we have already seen
(to account for Borel colorings).

Let us start with a few technical results about rigid surjections.
  For a chain $(A, \Boxed{<})$ and $x \in A$ let
  $\downarrow_A x = \{ y \in A: y \le x\}$. An \emph{initial segment} of $(A, \Boxed{<})$ is a subset
  $I \subseteq A$ such that $x \in I$ implies $\downarrow_A x \subseteq I$ for all $x \in A$.

  \begin{LEM}\label{rpemklei.lem.rs-basic-1}
    Let $(A, \Boxed{<})$ and $(B, \Boxed{<})$ be well ordered chains.
    A surjective map $f : A \to B$ is a rigid surjection from $(A, \Boxed{<})$ onto $(B, \Boxed{<})$
    if and only if $f$ takes every initial segment of $(A, \Boxed{<})$ onto an initial segment of $(B, \Boxed{<})$.
  \end{LEM}
  \begin{proof}
    Let us only show direction $(\Leftarrow)$. Take
    $b_1, b_2 \in B$ such that $b_1 < b_2$ and suppose that $\min f^{-1}(b_1) > \min f^{-1}(b_2) = a$.
    Let $I = \downarrow_A a$. Then $f(I) \ni b_2$, but $f(I) \not\ni b_1$ because $\min f^{-1}(b_1) > a$,
    whence $f^{-1}(b_1) \cap I = \0$. Therefore, $f(I)$ is not an initial segment of~$(B, \Boxed<)$.
  \end{proof}

  \begin{LEM}\label{rpemklei.lem.rs-basic-2}
    Let $g : A \to B$ and $h : B \to C$ be surjections. Assume that $(A, \Boxed{<})$, $(B, \Boxed{<})$ and $(C, \Boxed{<})$
    are well-ordered chains and that $g$ and $f = h \circ g$ are rigid surjections. Then $h$ is also
    a rigid surjection.
  \end{LEM}
  \begin{proof}
    Let us show that $h$ takes every initial segment of $(B, \Boxed{<})$ onto an initial segment of $(C, \Boxed{<})$.
    Let $I \subseteq B$ be an initial segment of $B$. Define $J \subseteq A$ as follows:
    $$
      J = \bigcup \{ \downarrow_A x : x \in A \text{ and } g(\downarrow_A x) \subseteq I \}.
    $$
    It is clear that $J$ is an initial segment of $(A, \Boxed{<})$ and that $g(J) \subseteq I$. To show that $g(J) = I$ take any $b \in I$
    and let $a = \min g^{-1}(b)$. To show that $a \in J$ it suffices to show that $g(\downarrow_A a) \subseteq I$. Take any $x \in \downarrow_A a$
    and let us show that $g(x) \le b$. Suppose this is not the case. Then $b < g(x)$. Since $g$ is rigid, we have that
    $a = \min g^{-1}(b) < \min g^{-1}(g(x)) \le x \le a$. Contradiction. This completes the proof that $g(J) = I$.
    Now, $h(I) = h(g(J)) = f(J)$, which is an initial segment of $(C, \Boxed{<})$ because $f$ is a rigid surjection.
  \end{proof}

  \begin{LEM}\label{rpemklei.lem.rs-basic-3}
    Let $f : A \to B$ be a surjection and let $(A, \Boxed{<})$ be a well-ordered chain. Define
    $f^\partial : B \to A$ by $f^\partial(b) = \min f^{-1}(b)$.

    $(a)$ Assume that $(B, \Boxed<)$ is well-ordered. Then
    $f$ is a rigid surjection if and only if $f^\partial : (B, \Boxed{<}) \to (A, \Boxed{<})$ is an embedding.
    
    $(b)$ There is a unique well-ordering of $B$ which turns $f$ into a rigid surjection.
  \end{LEM}
  \begin{proof} $(a)$ is obvious. Let us show $(b)$.
    According to $(a)$, in order to turn $f$ into a rigid surjection, $f^\partial$ has to be an embedding, and
    the unique linear order on $B$ that turns $f^\partial$ into an embedding (and hence $f$ into a rigid surjection) is given by
    $b_1 < b_2 \Leftrightarrow f^\partial(b_1) < f^\partial(b_2)$. The fact that $f^\partial$ is an embedding implies that
    $(B, \Boxed<)$ is isomorphic to a suborder of $(A, \Boxed<)$, so it has to be well-ordered.
  \end{proof}

\begin{LEM}\label{rpemklei.lem.lex-rs}
  Let $(A, \Boxed{<})$ and $(B, \Boxed{<})$ be well-ordered chains and let $f : A \to B$ be a rigid surjection
  $(A, \Boxed{<}) \to (B, \Boxed{<})$. Let $n$ be any positive integer and let $g : A^n \to B^n$ be the mapping
  defined by $g(a_1, \dots, a_n) = (f(a_1), \dots, f(a_n))$. Then:
  
  $(a)$ $\min g^{-1}(b_1, \dots, b_n) = (\min f^{-1}(b_1), \dots, \min f^{-1}(b_n))$, for all $b_1, \dots, b_n \in B$.
  
  $(b)$ $g : A^n \to B^n$ is a rigid surjection $(A^n, \Boxed{\lex{A}}) \to (B^n, \Boxed{\lex{B}})$.
\end{LEM}
\begin{proof}
  $(a)$
  Take any $b_1, \dots, b_n \in B$ and let $(a_1, \dots, a_n) = \min g^{-1}(b_1, \dots, b_n)$. We claim that
  $
    (a_1, \dots, a_n) = (\min f^{-1}(b_1), \dots, \min f^{-1}(b_n)).
  $
  Suppose this is not the case. Then there exists a $j \in \{1, 2, \dots, n\}$ such that
  $a_i = \min f^{-1}(b_i)$ for all $i < j$ and $a_j \ne \min f^{-1}(b_j)$. Note that, by construction, $f(a_j) = b_j$, that is,
  $a_j \in f^{-1}(b_j)$. Let $c = \min f^{-1}(b_j)$. Then $c < a_j$, so
  $$
    (a_1, \dots, a_{j-1}, c, a_{j+1}, \dots, a_n) \lex{A} (a_1, \dots, a_{j-1}, a_j, a_{j+1}, \dots, a_n)
  $$
  and
  $$
    (a_1, \dots, a_{j-1}, c, a_{j+1}, \dots, a_n) \in g^{-1}(b_1, \dots, b_n).
  $$
  But this contradicts the fact that $(a_1, \dots, a_n) = \min g^{-1}(b_1, \dots, b_n)$.
  
  $(b)$
  Take any $b_1, \dots, b_n, d_1, \dots, d_n \in B$ such that
  $
    (b_1, \dots, b_n) \lex{B} (d_1, \dots, d_n)
  $
  and let us show that
  $$
    \min g^{-1}(b_1, \dots, b_n) \lex{A} \min g^{-1}(d_1, \dots, d_n).
  $$
  By~$(a)$ it suffices to show that
  $$
    (\min f^{-1}(b_1), \dots, \min f^{-1}(b_n)) \lex{A} (\min f^{-1}(d_1), \dots, \min f^{-1}(d_n)).
  $$
  Since
  $
    (b_1, \dots, b_n) \lex{B} (d_1, \dots, d_n)
  $
  there is a $j \in \{1, \dots, n\}$ such that $b_i = d_i$ for $i < j$ and $b_j < d_j$. But then
  $\min f^{-1}(b_i) = \min f^{-1}(d_i)$ for $i < j$ and $\min f^{-1}(b_j) \lt{A} \min f^{-1}(d_j)$ because $f : A \to B$ is a rigid surjection.
  This completes the proof.
\end{proof}

\begin{CONSTR}\label{rpemkl.constr.1}
We shall now upgrade the functor $T$ and multiplication $\mu$
from Example~\ref{rpemklei.ex.T-Sigma} defined on $\Set^+$ to the category $\ChRs$ of well-ordered chains and rigid surjections.
Let $\Omega = \Omega_C \oplus \Omega_F$ be a well-ordered algebraic language where $\Omega_C$ is a well-ordered set of
constant symbols and $\Omega_F$ is a well-ordered set of functional symbols.
For a well-ordered chain $\calX = (X, \Boxed{<})$ of variables
let $T(X)$ be the set of all the $\Omega$-terms over the set of variables $X$.
We are now going to construct a particular well-ordering of $T(X)$ which is functorial and ensures that
every structure map of every $\Omega$-algebra is a rigid surjection.

Let $\xi = \min \calX$. For a term $t \in T(X)$, the \emph{shape of $t$} is the term obtained from
$t$ by replacing every variable that occurs in $t$ with $\xi$.
For example, for $\Omega = \{c, f, g\}$ where $f$ is ternary, $g$ is
binary and $c$ is a constant, the shape of $f(g(x_2, x_1), c, x_1)$ is $f(g(\xi, \xi), c, \xi)$,
while the shape of $x_1$ is~$\xi$. For each $n \ge 1$ let
$S_n$ be the set of all the shapes of length $n$ (recall that shapes are terms, and hence, strings of symbols).
Let us order $S_n$ lexicographically with respect to the well-ordering
$
  \{\xi\} \oplus \Omega \oplus \{\fbox{(} < \fbox{,} < \fbox{)}\},
$
where $\fbox{(}$, $\fbox{,}$ and $\fbox{)}$ are the usual additional symbols we use to
form terms. Note that each $S_n$ is well-ordered and that $S_1 = \{\xi\} \oplus \Omega_C$.

For every shape $\sigma \in \bigcup_{n \ge 1} S_n$ let $T_\sigma(X)$ denote the set of all the terms $t \in T(X)$ of
shape~$\sigma$. Let us order $T_\sigma(X)$ lexicographically with respect to the well-ordering
$\calX \oplus \Omega \oplus \{\fbox{(} < \fbox{,} < \fbox{)}\}$ and let us order $T(X)$ as follows:
$$
  (T(X), \Boxed<) = \left(\bigoplus_{\sigma \in S_1} T_\sigma(X) \right) \oplus \left(\bigoplus_{\sigma \in S_2} T_\sigma(X) \right) \oplus \left(\bigoplus_{\sigma \in S_3} T_\sigma(X) \right) \oplus \dots
$$
Note that this is a well-ordering of $T(X)$ and that
$$
  (T(X), \Boxed<) = \calX \oplus \Omega_C \oplus \left(\bigoplus_{\sigma \in S_2} T_\sigma(X) \right) \oplus \left(\bigoplus_{\sigma \in S_3} T_\sigma(X) \right) \oplus \dots
$$
We shall refer to this well-ordering of $T(X)$ as the \emph{a neat well-ordering of $T(X)$}.
Note that the neat well-ordering of $T(X)$ depends on $\calX$ and the choice of the well-ordering of~$\Omega$.

For a chain $\calX$ put $\hat T(\calX) = (T(X), \Boxed{<})$ where $<$ is the neat well-ordering of~$T(X)$.
This is how $\hat T$ acts on objects. For a rigid surjection $f : (X, \Boxed{<}) \to (Y, <)$ let
$\hat T(f) = T(f)$, that is, variable substitution where $\hat T(f)(t)$ is a new term
obtained from $t$ by systematically replacing each occurrence of $x \in X$ by $f(x) \in Y$. This will clearly be a functor
once we ensure that $\hat T(f)$ is well defined.

\begin{LEM}\label{rpemklei.lem.hat-T-rs}
  Let $(X, \Boxed{<})$ and $(Y, \Boxed{<})$ be well-ordered chains and let $f : (X, \Boxed{<}) \to (Y, \Boxed{<})$ be a rigid surjection.
  Then $\hat T(f) : \hat T(X) \to \hat T(Y)$ is a rigid surjection $(T(X), \Boxed{<}) \to
  (T(Y), \Boxed{<})$.
\end{LEM}
\begin{proof}
  Let us represent the chains $(T(X), \Boxed{<})$ and $(T(Y), \Boxed{<})$ as
  \begin{align*}
    (T(X), \Boxed{<}) &= (X, \Boxed{<}) \oplus \Omega_C \oplus \bigoplus_{n \ge 2}\left(\bigoplus_{\sigma \in S_n} T_\sigma(X) \right)\\
    (T(Y), \Boxed{<}) &= (Y, \Boxed{<}) \oplus \Omega_C \oplus \bigoplus_{n \ge 2}\left(\bigoplus_{\sigma \in S_n} T_\sigma(Y) \right).
  \end{align*}
  Being just a variable substitution, $\hat T(f)$ preserves the shape of terms, so $\hat T(f)$ restricted to $(T_\sigma(X), \Boxed{<})$
  for a fixed shape $\sigma$ is a rigid surjection $(T_\sigma(X), \Boxed{<}) \to (T_\sigma(Y), \Boxed{<})$.
  Note, next, that for each shape $\sigma$ the chain $(T_\sigma(X), \Boxed{<})$ is isomorphic to $(X^{n(\sigma)}, \Boxed{\lex{X}})$
  where $n(\sigma)$ is the number of occurrences of $\xi$ in $\sigma$. Lemma~\ref{rpemklei.lem.lex-rs} now implies that
  $\hat T(f)$ restricted to $(T_\sigma(X), \Boxed{<})$ is a rigid surjection for every shape $\sigma$, so
  $\hat T(f)$ as a whole is also a rigid surjection.
\end{proof}

Next, let us show that the multiplication $\hat \mu : \hat T \hat T \to \hat T$ defined as in Example~\ref{rpemklei.ex.T-Sigma}
remains well-defined. In other words, let us show that for every well-ordered chain $\calX = (X, \Boxed{<})$ the map
$\hat \mu_\calX : \hat T \hat T(\calX) \to \hat T(\calX)$ is a rigid surjection. To see that this is indeed the case recall that the
neat well-ordering of $\hat T(\calX)$ places all the variables before the terms involving at least one functional or constant symbol:
$$
  \hat T(\calX): \underbrace{x_1 < x_2 < \dots}_{\text{variables}} < \underbrace{t_1 < t_2 < \dots}_{\text{other terms}}
$$
In $\hat T \hat T(\calX)$ we use $\hat T(\calX)$ as variables upon which we build terms from $\hat T \hat T(\calX)$.
If we denote the elements of $\hat T(\calX)$ as
$$
  \langle x_1 \rangle < \langle x_2 \rangle < \dots < \langle t_1 \rangle < \langle t_2 \rangle < \dots
$$
just as a notational convenience, then the ordering of $\hat T \hat T(\calX)$ takes the following form:
$$
  \hat T \hat T(\calX): \underbrace{\langle x_1 \rangle < \langle x_2 \rangle < \dots < \langle t_1 \rangle < \langle t_2 \rangle < \dots}_{\text{variables}} < \underbrace{t'_1 < t'_2 < \dots}_{\text{other terms}}
$$
Recall also that $\mu_X$ acts as substitution of terms for variables, for example, $\mu_X(\langle t \rangle) = t$
and $\mu_X(t'(\langle t_1 \rangle, \langle t_2 \rangle)) = t'(t_1, t_2)$. Therefore,
the diagram depicting the action of $\hat \mu_\calX$ looks like this:
\begin{center}
  \begin{tikzcd}[column sep=-1.25mm]
    \hat T \hat T(\calX) \arrow[d,"\hat \mu_\calX"'] &:& \langle x_1 \rangle \arrow[d, mapsto] & \langle x_2 \rangle \arrow[d, mapsto]& \dots \arrow[d, mapsto] & \langle t_1 \rangle \arrow[d, mapsto]& \langle t_2 \rangle \arrow[d, mapsto]& \dots \arrow[d, mapsto] & t_1' \arrow[dlll, mapsto] & t_2' \arrow[dll, mapsto] & \dots \arrow[dllll, mapsto]\\
             \hat T(\calX)                    &:&         x_1         &         x_2         & \dots &         t_1         &         t_2         & \dots 
  \end{tikzcd}
\end{center}
which is clearly a rigid surjection.
\end{CONSTR}

A \emph{well-ordered $\Omega$-algebra} is a structure $\calA = (A, \Omega^A, \Boxed{<})$ where
$(A, \Omega^A)$ is an $\Omega$-algebra and $<$ is a well-ordering of $A$.
A mapping $f : A \to B$ is a \emph{rigid epimorphism} from
a well-ordered $\Omega$-algebra $\calA = (A, \Omega^A, \Boxed{<})$ onto a well-ordered
$\Omega$-algebra $\calB = (B, \Omega^B, \Boxed{<})$ if $f$ is a rigid surjection $(A, \Boxed{<}) \to (B, \Boxed{<})$
and at the same time an epimorphism $(A, \Omega^A) \to (B, \Omega^B)$.
Every well-ordered $\Omega$-algebra $\calA = (A, \Omega^A, \Boxed{<})$ is
also an Eilenberg-Moore $\hat T$-algebra $((A, \Boxed{<}), \eval^A)$ where the structure map
$$
  \eval^A : \hat T(A, \Boxed{<}) \to (A, \Boxed{<})
$$
is the evaluation in $\calA$, and vice versa.
Namely, it is easy to show that every structure map $\eval^A$ is a rigid surjection:
since the chain $\hat T(A, \Boxed{<})$ has the following form
$$
  \hat T(A, \Boxed{<}): \underbrace{a_1 < a_2 < \dots}_{\text{elements of $A$}} < \underbrace{t_1 < t_2 < \dots}_{\text{other terms}}
$$
and since $\eval^A$ is an evaluation of terms in $\calA$, it follows that the diagram depicting the action of $\eval^A$ looks like this:
\begin{center}
  \begin{tikzcd}[column sep=-1.25mm]
    \hat T(A, \Boxed{<}) \arrow[d,"\eval^A"'] &:&  a_1 \arrow[d, mapsto]& a_2 \arrow[d, mapsto]& \dots \arrow[d, mapsto]& t_1 \arrow[dlll, mapsto] & t_2 \arrow[dlll, mapsto] & \dots \arrow[dlll, mapsto]\\
             (A, \Boxed{<})                     &:&  a_1                  & a_2                  & \dots 
  \end{tikzcd}
\end{center}
which is clearly a rigid surjection.

  Since the category $\Alg(\Omega)$ of $\Omega$-algebras is isomorphic to the category of Eilenberg-Moore $T$-algebras (see Example~\ref{rpemklei.ex.T-Sigma})
  in what follows we shall think of an $\Omega$-algebra $\calA$ as an Eilenberg-Moore $T$-algebra $(A, \eval^A)$ and, thus, assume
  that $\Alg(\Omega) = \EM(T, \mu, \eta)$. Hence,
  
\begin{ASSUMP}\label{rpemkl.assumption.1}
  In the rest of the paper we treat every variety $\VV$ of $\Omega$-algebras as a subcategory of $\EM(T, \mu, \eta)$.
\end{ASSUMP}

Let $\WalgRe(\Omega)$ denote the category of well-ordered $\Omega$-algebras
and rigid epimorphisms (that is, algebraic homomorphisms that are at the same time
rigid surjections) understood as a full subcategory of $\EM^w_\re(\hat T, \hat \mu)$.
For a variety $\VV$ of $\Omega$-algebras let $\WalgRe(\VV)$ denote the full subcategory of $\WalgRe(\Omega)$
spanned by all the well-ordered $\Omega$-algebras $\calA = ((A, \Boxed{<}), \eval^A)$ where $(A, \eval^A) \in \Ob(\VV)$.
Hence, for every variety $\VV$ of $\Omega$-algebras
$\WalgReFin(\VV)$ is isomorphic to a subcategory of the category $\EM^w_\re(\hat T, \hat \mu)$.

Let $\VV$ be a nontrivial variety of $\Omega$-algebras and $\calX = (X, \Boxed{<})$ a well-ordered chain.
Recall that $\nu_{\VV,X} : \calF(X) \to \calF_\VV(X)$ is a natural epimorphism of the term algebra $\calF(X) = (T(X), \mu_X)$
onto the free $\VV$ algebra $\calF_\VV(X) = (T_\VV(X), \theta_{\VV, X})$, where $\theta_{\VV, X}$ denotes the structure map of
$\calF_\VV(X)$.
By Lemma~\ref{rpemklei.lem.rs-basic-3} there is a unique well-ordering of $T_\VV(X)$ which turns $\nu_{\VV,X}$ into a rigid surjection
$\nu_{\VV,X} : (T(X), \Boxed{<}) \to (T_\VV(X), \Boxed{<})$. Let $\hat T_\VV(\calX)$ denote the chain $(T_\VV(X), <)$. Then
$\nu_{\VV,\calX} : (\hat T(\calX), \hat \mu_\calX) \to (\hat T_\VV(\calX), \theta_{\VV, X})$ is a rigid epimorphism.
Let $\hat \calF(\calX) = (\hat T(\calX), \hat \mu_\calX)$ and
$\hat \calF_\VV(\calX) = (\hat T_\VV(\calX), \theta_{\VV, X})$.

\begin{LEM}\label{rpemklei.lem.f-star}
  Let $\VV$ be a nontrivial variety of $\Omega$-algebras and let $\calX = (X, \Boxed<)$ and $\calY = (Y, \Boxed<)$
  be well-ordered chains. For every $f \in \hom_{\ChRs}(\calX, \calY)$ there is an
  $f^* \in \hom_{\WalgRe(\VV)}(\hat \calF_\VV(\calX), \hat \calF_\VV(\calY))$ such that the diagram below commutes in
  $\WalgRe(\Omega)$:
  \begin{center}
    \begin{tikzcd}
      \hat \calF(\calX) \ar[r, "\hat T(f)"] \ar[d, "\nu_{\VV, \calX}"'] & \hat \calF(\calY) \ar[d, "\nu_{\VV, \calY}"] \\
      \hat \calF_\VV(\calX) \ar[r, "f^*"] & \hat \calF_\VV(\calY)
    \end{tikzcd}
  \end{center}
  We shall denote $f^*$ by $\hat T_\VV(f)$.
\end{LEM}
\begin{proof}
  In the unordered setting, every map $f : X \to Y$ between the generating sets determines homomorphisms
  $T(f) : \calF(X) \to \calF(Y)$ and $f^* : \calF_\VV(X) \to \calF_\VV(Y)$ between the corresponding free algebras.
  Both $T(f)$ and $f^*$ are epimorphisms and the following diagram commutes:
  \begin{center}
    \begin{tikzcd}
      \calF(X) \ar[r, "T(f)"] \ar[d, "\nu_{\VV, X}"'] & \calF(Y) \ar[d, "\nu_{\VV, Y}"] \\
      \calF_\VV(X) \ar[r, "f^*"] & \calF_\VV(Y)
    \end{tikzcd}
  \end{center}
  Moving back to the ordered setting, $\hat T(f)$ is a rigid surjection by Lemma~\ref{rpemklei.lem.hat-T-rs}.
  Since $\nu_{\VV, \calX}$ and $\nu_{\VV, \calY}$ are rigid surjections by construction, $f^*$ is a rigid surjection
  by Lemma~\ref{rpemklei.lem.rs-basic-2}.
\end{proof}

The following lemma confirms the intuition that the natural surjection $\nu_{\VV,\omega}$ remains a reflection
even in case of ordered algebras. This is an important technical prerequisite for establishing the dual Ramsey property
in lemmas that follow. Recall that (co)reflections were defined immediately after Lemma~\ref{rpemklei.lem.C-D}.

\begin{LEM}\label{rpemklei.lem.TAlg-nu}
  Let $\VV$ be a nontrivial variety of $\Omega$-algebras. Then
  $$
    \nu_{\VV,\omega} : (\hat T(\omega), \hat \mu_{\omega}) \to (\hat T_\VV(\omega), \theta_{\VV, \omega})
  $$
  is a reflection of $(\hat T(\omega), \hat \mu_{\omega})$ in $\WalgRe(\VV)$.
\end{LEM}
\begin{proof}
  Let $((A, \Boxed{<}), \alpha)$ be a well-ordered algebra such that $(A, \alpha) \in \VV$
  and let $f : (\hat T(\omega), \hat\mu_\omega) \to ((A, \Boxed{<}), \alpha)$ be a rigid epimorphism.  
  It is a well-known fact that in the usual (unordered) setting
  $\nu_{\VV,\omega} : (T(\omega), \mu_\omega) \to (T_\VV(\omega), \theta_{\VV,\omega})$
  is a reflection of $(T(\omega), \mu_\omega)$ in $\VEpi$.
  Therefore, there is an epimorphism $g : (T_\VV(\omega), \theta_{\VV, \omega}) \to (A, \alpha)$ such that the following diagram
  commutes in $\AlgEpi(\Omega)$:
  \begin{center}
    \begin{tikzcd}
             &  (T(\omega), \mu_\omega) \arrow[d, "\nu_{\VV,\omega}"] \arrow[dl, "f"']\\
        (A, \alpha) & (T_\VV(\omega), \theta_{\VV,\omega}) \arrow[l, "g"]
    \end{tikzcd}
  \end{center}
  Going back to the well-ordered setting, $f : (\hat T(\omega), \hat \mu_\omega) \to ((A, \Boxed{<}), \alpha)$
  and $\nu_{\VV,\omega} : (\hat T(\omega), \hat\mu_\omega) \to (\hat T_\VV(\omega), \theta_{\VV,\omega})$
  are rigid surjections, so Lemma~\ref{rpemklei.lem.rs-basic-2}
  ensures that $g : (T_\VV(\omega), \theta_{\VV,\omega}) \to (A, \alpha)$ is not only an epimorphism, but also a rigid surjection
  $(\hat T_\VV(\omega), \theta_{\VV,\omega}) \to ((A, \Boxed{<}), \alpha)$. Hence, $g$ is a rigid epimorphism.
\end{proof}

\begin{LEM}\label{rpemklei.lem.TAlg-1}
  Let $\VV$ be a nontrivial variety of\/ $\Omega$-algebras.
  Taking $\EM = \EM^w_\re(\hat T, \hat \mu)$ as the ambient category,
  for all $\calA, \calB \in \WalgReFin(\VV)$ and all $k \ge 2$ we have that $(\hat T_\VV(\omega), \theta_{\VV,\omega}) \longleftarrow (\calB)^\calA_k$.
\end{LEM}
\begin{proof}
  Take any $\calA, \calB \in \WalgReFin(\VV)$ and any $k \ge 2$, and let us first show that
  $(\hat T(\omega), \hat\mu_\omega) \longleftarrow (\calB)^\calA_k$.
  To do so let us construct a pre-adjunction
  $$
    F : \Ob(\EM^\op) \rightleftarrows \Ob(\ChRs^\op) : H.
  $$
  For $\calB = ((B, \Boxed{<}), \beta) \in \Ob(\EM)$ put $F(\calB) = (B, \Boxed{<})$,
  for $(C, \Boxed{<}) \in \Ob(\ChRs)$ put $H(C, \Boxed{<}) = (\hat T(C, \Boxed{<}), \hat\mu_{(C, \Boxed{<})})$ and
  for $u \in \hom_{\ChRs}((C, \Boxed{<}), (B, \Boxed{<}))$ put $\Phi_{\calB, (C, \Boxed{<})}(u) = \beta \cdot \hat T(u)$.
  By dualizing the proof of Theorem~\ref{rpemklei.thm.ramsey} we can now easily show that
  $F : \Ob(\EM^\op) \rightleftarrows \Ob(\ChRs^\op) : H$ is indeed a pre-adjunction. 
  
  Since $\ChRsFin$ has the dual Ramsey property
  there is a finite chain $(C, \Boxed{<})$ such that $(C, \Boxed{<}) \longleftarrow (F(\calB))^{F(\calA)}_k$.
  By Theorem~\ref{opos.thm.main}~$(a)$ it then follows that $H(C, \Boxed{<}) = (\hat T(C, \Boxed{<}), \hat\mu_{(C, \Boxed{<})}) \longleftarrow (\calB)^{\calA}_k$.
  Now, take any rigid surjection $f : (\omega, \Boxed<) \to (C, \Boxed{<})$. The fact that $\hat \mu$ is natural
  yields that $\hat T(f) : \hat T(\omega) \to \hat T(C, \Boxed{<})$
  is a morphism in $\EM$. The dual of Lemma~\ref{rpemklei.lem.C-D} now
  ensures that $(\hat T(\omega), \hat\mu_\omega) \longleftarrow (\calB)^{\calA}_k$.

  We have shown in Lemma~\ref{rpemklei.lem.TAlg-nu}
  that $\nu_{\VV,\omega} : (\hat T(\omega), \hat\mu_\omega) \to (\hat T_\VV(\omega), \theta_{\VV,\omega})$
  is a reflection of $(\hat T(\omega), \hat\mu_\omega)$ in $\WalgRe(\VV)$.
  Hence, the dual of Lemma~\ref{rpemklei.lem.refl} now yields that $(\hat T_\VV(\omega), \theta_{\VV,\omega}) \longleftarrow (\calB)^\calA_k$.
\end{proof}

As in the proof of Theorem~\ref{rpemklei.thm.ramsey}, the pre-adjunction constructed here
cannot be expected to be an adjunction. Although $F$ and $H$ behave as the forgetful functor and the free-algebra functor, respectively,
one cannot hope for $\Phi$ to be either bijective on hom-sets or natural in its arguments.

\begin{LEM}\label{rpemklei.lem.TAlg-2}
  Let $\VV$ be a nontrivial variety of $\Omega$-algebras.
  Taking $\EM = \EM^w_\re(\hat T, \hat \mu)$ as the ambient category, $(\hat T_\VV(\omega), \theta_{\VV,\omega})$
  is projectively universal and projectively weakly locally finite for $\WalgReFin(\VV)$.
\end{LEM}
\begin{proof}
  To see that $(\hat T_\VV(\omega), \theta_{\VV,\omega})$ is projectively universal for $\WalgReFin(\VV)$
  take any $\calB = ((B, \Boxed{<}), \beta) \in \WalgReFin(\VV)$ and any rigid surjection $f : \omega \to (B, \Boxed{<})$.
  Then the square on the left commutes because $\hat \mu$ is natural,
  while the square on the right commutes because $\calB$ is a weak Eilenberg-Moore $\hat T$-algebra:
  \begin{center}
    \begin{tikzcd}
          \hat T \hat T(\omega)  \arrow[d, "\hat \mu_{\omega}"']   \arrow[rr, "\hat T \hat T(f)"]
      & & \hat T \hat T(B, \Boxed{<}) \arrow[d, "\hat \mu_{(B, \Boxed<)}"']  \arrow[rr, "\hat T(\beta)"]
      & & \hat T(B, \Boxed{<})          \arrow[d, "\beta"]
    \\
          \hat T(\omega)    \arrow[rr, "\hat T(f)"']
      & & \hat T(B, \Boxed{<})   \arrow[rr, "\beta"']
      & & (B, \Boxed{<})
    \end{tikzcd}
  \end{center}
  Therefore, $(\hat T(\omega), \hat\mu_\omega) \overset\EM\longrightarrow \calB$.
  Now recall that $\nu_{\VV,\omega} : (\hat T(\omega), \hat\mu_\omega) \to (\hat T_\VV(\omega), \theta_{\VV,\omega})$
  is a reflection of $(\hat T(\omega), \hat\mu_\omega)$
  in $\WalgRe(\VV)$ (Lemma~\ref{rpemklei.lem.TAlg-nu}). This ensures the
  existence of a morphism $(\hat T_\VV(\omega), \theta_{\VV,\omega}) \overset\EM\longrightarrow \calB$.

  To see that $(\hat T_\VV(\omega), \theta_{\VV,\omega})$ is projectively weakly locally finite for $\WalgReFin(\VV)$ take any
  $\calA = ((A, \Boxed{<}), \alpha)$ and $\calB = ((B, \Boxed{<}), \beta)$ in $\WalgReFin(\VV)$ and
  arbitrary morphisms $f : (\hat T_\VV(\omega), \theta_{\VV,\omega}) \to \calA$ and 
  $g : (\hat T_\VV(\omega), \theta_{\VV,\omega}) \to \calB$. Then
  $f : (T_\VV(\omega), \theta_{\VV,\omega}) \to (A, \alpha)$ and $g : (T_\VV(\omega), \theta_{\VV,\omega}) \to (B, \beta)$
  are epimorphisms in the category $\AlgEpi(\Omega)$ of (unordered) $\Omega$-algebras.
  The category $\AlgEpi(\Omega)$ has products, and if the two algebras come from a variety $\VV$ the product also belongs to $\VV$.
  So, let $(A \times B, \delta) = (A, \alpha) \times (B, \beta)$ be the
  product of the two finite algebras. Clearly, there is a unique homomorphism $h_0 : (T_\VV(\omega), \theta_{\VV,\omega}) \to (A \times B, \delta)$
  such that $\pi_1 \circ h_0 = f$ and $\pi_2 \circ h_0 = g$, but $h_0$ is not necessarily an epimorphism.
  However, $C = h_0(T(\omega))$ is a subalgebra of $(A \times B, \delta)$. Let $\gamma : T(C) \to C$ be the corresponding
  structure map. We now have that the codomain restriction $h : (T_\VV(\omega), \theta_{\VV,\omega}) \to (C, \gamma)$ of $h_0$ is an
  epimorphism and the following diagram commutes in $\AlgEpi(\Omega)$:
  \begin{center}
    \begin{tikzcd}
      (T_\VV(\omega), \theta_{\VV,\omega}) \arrow[rr, "h"] \arrow[d, "f"'] \arrow[drr, "g"', near end] & & (C, \gamma) \arrow[dll, "\pi_1", near end] \arrow[d, "\pi_2"] \\
      (A, \alpha) & & (B, \beta)
    \end{tikzcd}
  \end{center}
  Note that $(C, \gamma)$ is finite as a subalgebra of a finite algebra.
  According to Lemma~\ref{rpemklei.lem.rs-basic-3} there is a unique well-ordering of $C$ which turns $h$ into a rigid surjection,
  and hence into a rigid epimorphism  $h : (\hat T_\VV(\omega), \theta_{\VV,\omega}) \to ((C, \Boxed{<}), \gamma)$.
  Lemma~\ref{rpemklei.lem.rs-basic-2} now ensures that $\pi_1 : ((C, \Boxed{<}), \gamma) \to \calA$ and
  $\pi_2 : ((C, \Boxed{<}), \gamma) \to \calB$ are also rigid surjections, and hence rigid epimorphisms.
  Therefore, $(\hat T_\VV(\omega), \theta_{\VV,\omega})$ is projectively locally finite for $\WalgReFin(\VV)$.
\end{proof}

\begin{THM}\label{rpemklei.thm.ALG-MAIN}
  Let $\Omega$ be an arbitrary algebraic language and $\VV$ a nontrivial variety of $\Omega$-algebras.
  Let $\CC$ be the category whose objects are ordered finite algebras $(A, \Omega^A, \Boxed<)$
  where $(A, \Omega^A) \in \VV^\fin$ and $<$ is a linear ordering of $A$,
  and whose morphisms are rigid epimorphisms (that is, epimorphisms between algebras that are at the
  same time rigid surjections). Then $\CC$ has the dual Ramsey property.
\end{THM}
\begin{proof}
  Having in mind Construction~\ref{rpemkl.constr.1} and Assumption~\ref{rpemkl.assumption.1},
  without loss of generality we may take that $\CC = \WalgReFin(\VV)$.
  Take $\EM = \EM^w_\re(\hat T, \hat \mu)$ as the ambient category. We have seen in Lemma~\ref{rpemklei.lem.TAlg-2} that
  $(\hat T_\VV(\omega), \theta_{\VV,\omega})$ is projectively universal and projectively locally finite for $\CC$.
  Lemma~\ref{rpemklei.lem.TAlg-1} shows that
  for all $\calA, \calB \in \Ob(\CC)$ and all $k \ge 2$ we have that $(\hat T_\VV(\omega), \theta_{\VV,\omega}) \longleftarrow (\calB)^\calA_k$.
  Therefore, by the dual of Lemma~\ref{akpt.lem.ramseyF} we have that $t^\partial_{\CC}(\calA) = 1$ for all $\calA \in \Ob(\CC)$.
  This is just another way of saying that $\CC$ has the dual Ramsey property.
\end{proof}

\begin{COR}\label{rpemklei.cor.ALG-MAIN}
  For every algebraic language $\Omega$ and every nontrivial variety $\VV$ of $\Omega$-algebras
  the category $\VV^\fin_\epi$ of finite $\VV$ algebras and epimorphisms has finite dual small Ramsey degrees.
\end{COR}
\begin{proof}
  As a matter of notational convenience put $\CC = \WalgReFin(\VV)$ and $\BB = \VEpiFin$.
  Let $U : \CC^\op \to \BB^\op$ be the forgetful functor that forgets the order and let us
  show that $U$ is a reasonable expansion with unique restrictions.
  
  To see that $U$ is reasonable let $(A, \alpha), (B, \beta) \in \Ob(\BB)$ be
  finite $\VV$ algebras, let $e : (B, \beta) \to (A, \alpha)$
  be an epimorphism and let $<$ be an arbitrary linear order on $A$:
  \begin{center}
    \begin{tikzcd}
      (A, \alpha, \Boxed{<}) \arrow[d, dashed, mapsto, "U"'] & & & &\CC^\op\\
      (A, \alpha) \arrow[rr, "e"'] & & (B, \beta) & &\BB^\op
    \end{tikzcd}
  \end{center}
  (note that an epimorphism $e : B \to A$ in $\BB$ is an arrow $A \to B$ in $\BB^\op$).
  Then it is easy to find a linear order on $B$ such that $e : (B, \Boxed{<}) \to (A, \Boxed{<})$
  is a rigid surjection. This turns $e$ into a rigid epimorphism $e : (B, \beta, \Boxed{<}) \to (A, \alpha, \Boxed{<})$.
  
  Let us now show that $U$ has unique restrictions.
  Let $(A, \alpha), (B, \beta) \in \Ob(\BB)$, let $e : (B, \beta) \to (A, \alpha)$
  be an epimorphism and let $<$ be an arbitrary linear order on $B$:
  \begin{center}
    \begin{tikzcd}
      & & (B, \beta, \Boxed{<}) \arrow[d, dashed, mapsto, "U"]   & &\CC^\op\\
      (A, \alpha) \arrow[rr, "e"'] & & (B, \beta) & &\BB^\op
    \end{tikzcd}
  \end{center}
  By Lemma~\ref{rpemklei.lem.rs-basic-3}
  there is a unique linear order $<$ on $A$ such that $e : (B, \Boxed{<}) \to (A, \Boxed{<})$ is a
  rigid surjection.
  
  Since $\CC^\op$ has the Ramsey property and $U^{-1}(A, \alpha)$ is finite for every $(A, \alpha)$
  (because there are only finitely many linear orders on a finite set) Theorem~\ref{sbrd.thm.small1} implies that
  $\BB^\op$ has finite small Ramsey degrees, so $\BB$ has finite dual small Ramsey degrees.
\end{proof}

So, although we still do not know whether the class of all finite groups has a precompact Ramsey expansion,
the following is a straightforward consequence of Theorem~\ref{rpemklei.thm.ALG-MAIN} and
Corollary~\ref{rpemklei.cor.ALG-MAIN}:

\begin{COR}\label{rpemklei.cor.groups}
  Let let $\VV$ be a nontrivial variety of groups.
  Then the category whose objects are ordered finite groups from $\VV$ and morphisms are rigid epimorphisms
  has the dual Ramsey property.
  Moreover, the category $\VEpiFin$ has finite dual small Ramsey degrees.
\end{COR}
  
Of course, the same is true for any variety of rings, modules,
lattices (in particular, distributive lattices, modular lattices, \dots),
boolean algebras, vector spaces and so on.

Finally, let us prove that for a countable algebraic language $\Omega$, in any variety $\VV$ of $\Omega$-algebras
finite algebras have finite dual big Ramsey degrees in the free $\VV$ algebra on $\omega$ generators.
As usual, we shall first prove the ordered version of the statement and then infer the
unordered version.

\begin{THM}\label{rpemklei.thm.main-V-order}
  Let $\Omega$ be a countable algebraic language and $\VV$ a nontrivial variety of $\Omega$-algebras.
  Let $\CC$ be the category whose objects are well-ordered algebras $(A, \Omega^A, \Boxed<)$
  where $(A, \Omega^A) \in \VV$ and $<$ is a well-ordering of $A$,
  and whose morphisms are rigid epimorphisms. Then 
  every finite ordered $\VV$ algebra has finite big dual Ramsey degree in $\hat\calF_\VV(\omega)$
  with respect to Borel colorings. More precisely, for every $\calA = (A, \Omega^A, \Boxed<) \in \CC^\fin$:
  $$
    \flatdualT_\CC(\calA, \hat\calF_\VV(\omega)) \le |A|.
  $$
\end{THM}
\begin{proof}
  Having in mind Construction~\ref{rpemkl.constr.1} and Assumption~\ref{rpemkl.assumption.1},
  without loss of generality we may take that $\CC = \WalgRe(\VV)$.
  Since $\VV$ is a nontrivial variety we have that $(i, j) \notin \Theta_\VV(\omega)$ whenever $i \ne j$
  (in other words, $\Theta_\VV(\omega)$ cannot identify distinct generators in $T(\omega)$ for otherwise $\VV$ would be a trivial variety),
  so the chain $\hat T_\VV(\omega)$ takes the following form:
  \begin{equation}\label{rpemkl.eq.hatTVomega}
     \hat T_\VV(\omega) = \omega \oplus \Omega_C \oplus \dots
  \end{equation}
  (actually, $\hat T_\VV(\omega)$ is isomorphic to the right-hand side because the
  elements of $\hat T_\VV(\omega)$ are congruence classes, but there is no harm in identifying
  the congruence class of a generator $x/\Theta_\VV(\omega) = \{x\}$ with the generator $x$ itself, $x \in \omega$).

  Take any $\calA = (A, \Omega^A, \Boxed<) \in \CC^\fin$ where
  $(A, \Boxed<) = \{a_1 < a_2 < \ldots < a_s\}$, $s = |A|$. As a notational convenience
  let
  $$
    F = T_\VV(\omega)
  $$
  be the carrier of the free algebra $\calF_\VV(\omega)$ (without the ordering relation), let
  $$
    H = \hom_{\Alg(\Omega)}(\calF_\VV(\omega), (A, \Omega^A)) \subseteq A^F
  $$
  be the set of all the homomorphisms $\calF_\VV(\omega) \to (A, \Omega^A)$, and let
  $$
    R = \hom_\CC(\hat\calF_\VV(\omega), \calA) \subseteq H
  $$
  be the set of all the rigid epimorphisms $\hat\calF_\VV(\omega) \to \calA$.
  Note that $F$ is a countable set and that $H$ is closed in $A^F$. Next, define
  $$
    \pi : H \to A^\omega \text{\quad by\quad} \pi(h) = \restr h \omega.
  $$
  
  \medskip

  Claim 1. $\pi$ is a continuous bijection.
  
  Proof. Since $\calF_\VV(\omega)$ is a free $\VV$ algebra, every mapping $g : \omega \to A$
  uniquely determines a homomorphism $g^\# : \calF_\VV(\omega) \to (A, \Omega^A)$.
  Therefore, $\pi$ is a bijection. It is clear that $\pi$ is continuous.
  
  \medskip
  
  Claim 2. Both $R$ and $H \setminus R$ are Borel in $A^F$.
  
  Proof. It suffices to show that $R$ is Borel in $H$.
  Let us denote the basic open sets in $A^F$ by
  $
    B\left(\begin{smallmatrix} x_1 & \dots & x_n \\ a_1 & \dots & a_n \end{smallmatrix}\right) =
    \{f \in A^F : f(x_i) = a_i, 1 \le i \le n\}
  $.
  It is easy to see that $x = \min f^{-1}(a)$ iff
  $$
    f \in \rho(x, a) =
    B\left(\begin{smallmatrix} x \\ a \end{smallmatrix}\right)
    \cap
    \bigcap_{\begin{smallmatrix}x' \in F\\x' < x\end{smallmatrix}} (A^F \setminus B\left(\begin{smallmatrix} x' \\ a \end{smallmatrix}\right)).
  $$
  Then, for $f \in H$ we have that $f \in R$ iff $f$ is a surjection, and a rigid one:
  $$
    f \in \Bigg(\bigcap_{a \in A} \bigcup_{x \in F} B\left(\begin{smallmatrix} x \\ a \end{smallmatrix}\right)\Bigg)
    \cap
    \Bigg(\bigcap_{a \in A} \bigcap_{\begin{smallmatrix}b \in A\\a < b\end{smallmatrix}} \bigcup_{x \in F} \bigcup_{\begin{smallmatrix}y \in F\\x < y\end{smallmatrix}} (\rho(x, a) \cap \rho(y, b))\Bigg).
  $$
  This concludes the proof of Claim 2.
  
  \medskip
  
  Claim 3. $\pi(R)$ is Borel in $A^\omega$. If $\{R_1, \ldots, R_k\}$ is a partition of $R$ into Borel sets, then
  $\{\pi(R_1), \ldots, \pi(R_k)\}$ is a partition of $\pi(R)$ into Borel sets.
  
  Proof. Follows from the Lusin-Suslin Theorem, having in mind that $\pi$ is a continuous bijection.
  
  \medskip

  Finally, for $1 \le i \le s$ let
  $$
    A_i = \{a_1 < \dots < a_i\} \text{\quad and\quad} S_i = \hom_{\ChRs}(\omega, A_i) \subseteq A^\omega.
  $$
  By the argument we used in the proof of Claim 2 it easily follows that each $S_i$ is Borel in $A^\omega$.
  The purpose of the sets $S_i$ is to make it easier for us to control the rigid surjections
  $\hat\calF_\VV(\omega) \to \calA$. Namely, the choice of the ordering of the carrier set of $\hat\calF_\VV(\omega)$
  as given in \eqref{rpemkl.eq.hatTVomega} suggests that every rigid surjection $\hat\calF_\VV(\omega) \to \calA$
  is uniquely determined by a rigid surjection of $\omega$ into an initial segment of~$A$. All such
  rigid surjections are, therefore, partitioned as $S_1 \cup \ldots \cup S_s$, and the proof then proceeds
  ``by parts''.
  
  \medskip
  
  We are now ready to show that $\flatdualT_\CC(\calA, \hat\calF_\VV(\omega)) \le |A|$. Take any $k \ge 2$
  and an arbitrary Borel coloring $\chi : R \to k$. Define $\gamma : \pi(R) \to k$ by
  $$
    \gamma(f) = \chi(\pi^{-1}(f)).
  $$
  Claim 3 now ensures that $\gamma$ is a Borel coloring. Finally,
  define Borel colorings $\gamma_i : S_i \to k$, $1 \le i \le s$, by
  $$
    \gamma_i(f) = \begin{cases}
      \gamma(f), & f \in \pi(R) \cap S_i,\\
      0, & \text{otherwise}.
    \end{cases}
  $$
  Let us construct
  $\gamma'_i : S_i \to k$ and $w_i \in \hom_{\ChRs}(\omega, \omega)$, $i \in \{1, \ldots, s\}$,
  inductively as follows. First, put $\gamma'_s = \gamma_s$.
  Given a Borel coloring $\gamma'_i : S_i \to k$, construct $w_i$
  by the Infinite Dual Ramsey Theorem (Theorem~\ref{rpemklei.thm.IDRT}):
  since $\omega \Borel\longleftarrow (\omega)^{A_i}_{k}$, there is a $w_i \in \hom_{\ChRs}(\omega, \omega)$
  such that
  $$
    |\gamma'_i(S_i \cdot w_i)| \le 1.
  $$
  %(We write the inequality here because we cannot guarantee that $S_i \ne \0$.)
  Finally, given $w_i \in \hom_{\ChRs}(\omega, \omega)$ define $\gamma'_{i-1} : S_{i-1} \to k$ by
  $$
    \gamma'_{i-1}(f) = \gamma_{i-1}(f \cdot w_i \cdot \dots \cdot w_s).
  $$
  Since $\gamma_{i-1}$ is a Borel coloring and the composition of morphisms is continuous, $\gamma'_{i-1}$ is also a Borel coloring.
  
  Now, put $u = w_1 \cdot \ldots \cdot w_s \in \hom_{\ChRs}(\omega,\omega)$ and let us show that
  $$
    |\chi(R \cdot \hat T_\VV(u))| \le s.
  $$
  (see Lemma~\ref{rpemklei.lem.f-star}). Since $\pi$ is a bijection,
  $$
    \chi(R \cdot \hat T_\VV(u)) = \chi(\pi^{-1}(\pi(R \cdot \hat T_\VV(u)))) = \gamma(\pi(R \cdot \hat T_\VV(u))).
  $$

  \medskip
  
  Claim 4. $\pi(R \cdot \hat T_\VV(u)) = \pi(R) \cdot u \subseteq \pi(R)$.
  
  Proof. Since $R \cdot \hat T_\VV(u) \subseteq R$ we immediately have that $\pi(R \cdot \hat T_\VV(u)) \subseteq \pi(R)$.
  To see that $\pi(R \cdot \hat T_\VV(u)) = \pi(R) \cdot u$ take any $h \in R$.
  By the construction of $T_\VV(u)$ (see Lemma~\ref{rpemklei.lem.f-star}) it follows that
  $\restr{(h \cdot \hat T_\VV(u))}{\omega} = \restr h \omega \cdot u$. This concludes the proof of the claim.
  
  \medskip
  
  Therefore, $\chi(R \cdot \hat T_\VV(u)) = \gamma(\pi(R) \cdot u)$. Since $\pi(R) \subseteq \bigcup_{i=1}^s S_i$,
  $$
    |\gamma(\pi(R) \cdot u)|
    \le |\gamma(\bigcup_{i=1}^s S_i \cdot u)|
    = |\bigcup_{i=1}^s \gamma(S_i \cdot u)|
    \le \sum_{i=1}^s |\gamma(S_i \cdot u)|.
  $$
  Fix an $i \in \{1, \dots, s\}$. Clearly, $S_i \cdot u \subseteq S_i$
  and $S_i \cdot w_1 \cdot \ldots \cdot w_i \subseteq S_i \cdot w_i$, whence
  $$
    |\gamma(S_i \cdot u)|
    = |\gamma_i(S_i \cdot w_1 \cdot \ldots \cdot w_s)|
    = |\gamma'_i(S_i \cdot w_1 \cdot \ldots \cdot w_i)|
    \le |\gamma'_i(S_i \cdot w_i)| \le 1.
  $$
  Putting it all together, we finally get
  $$
    |\chi(R \cdot \hat T_\VV(u))| = |\gamma(\pi(R) \cdot u)| \le \sum_{i=1}^s |\gamma(S_i \cdot u)| \le s.
  $$
  This completes the proof.
\end{proof}

\begin{COR}\label{rpemklei.cor.dual-big-rd-unordered}
  Let $\Omega$ be a countable algebraic language and $\VV$ a nontrivial variety of $\Omega$-algebras.
  Then every finite $\VV$ algebra has a finite big dual Ramsey degree
  with respect to Borel colorings in $\calF_\VV(\omega)$, the free
  $\VV$ algebra on $\omega$ generators, taking the category $\VEpi$ of $\VV$ algebras and epimorphisms as the ambient category.
  More precisely, for every $\calA \in \VEpiFin$ with $n$ elements:
  $$
    \flatdualT_{\VEpi}(\calA, \calF_\VV(\omega)) \le n \cdot n!
  $$
\end{COR}
\begin{proof}
  As a matter of notational convenience put $\CC = \WalgRe(\VV)$ and $\BB = \VEpi$.
  Let $U : \CC^\op \to \BB^\op$ be the forgetful functor that forgets the order.
  The argument used in Claim 2 of the proof of Theorem~\ref{rpemklei.thm.main-V-order}
  can be used here as well to show that $\hom_{\CC^\op}(\hat\calA, \hat\calB)$ is a Borel set in $\hom_{\BB^\op}(\calA, \calB)$.
  (Recall that we work with categories enriched over $\Top$.)
  To see that $U$ has unique restrictions we can repeat the argument from the proof of Corollary~\ref{rpemklei.cor.ALG-MAIN}.
  
  Let $\calA = (A, \Omega^A) \in \Ob(\BB^\fin)$ be an arbitrary finite $\VV$ algebra.
  As we have seen in Theorem~\ref{rpemklei.thm.main-V-order}, there exists a linear ordering $\hat\calF_\VV(\omega)$ of
  $\calF_\VV(\omega)$ such that every finite ordered $\VV$ algebra has a big dual Ramsey degree in $\hat\calF_\VV(\omega)$
  with respect to Borel colorings which does not exceed the number of elements of the algebra.
  Therefore, for every linear ordering $<$ of $A$ we have that
  $$
    \flatdualT_\CC((A, \Omega^A, <), \hat\calF_\VV(\omega)) \le |A|.
  $$
  Since $U^{-1}(\calA)$ is finite (because there are only finitely many linear orders on a finite set)
  Theorem~\ref{sbrd.thm.big1-dual-borel} tells us that
  $$
    \flatT_{\BB^\op}(\calA, \calF_\VV(\omega)) \le \sum_{\calA^* \in U^{-1}(\calA)} \flatT_{\CC^\op}(\calA^*, \hat\calF_\VV(\omega)).
  $$
  In other words,
  $$
    \flatdualT_{\BB}(\calA, \calF_\VV(\omega)) \le \sum_{\calA^* \in U^{-1}(\calA)} \flatdualT_{\CC}(\calA^*, \hat\calF_\VV(\omega)) \le n \cdot n!,
  $$
  where $n = |A|$. This completes the proof.
\end{proof}

\section{Acknowledgements}

The author would like to thank his colleagues Roz\'alia Madar\'asz and Milo\v s Kurili\'c for many helpful suggestions.

The author gratefully acknowledges the financial support of the Ministry of Education, Science and Technological Development
of the Republic of Serbia (Grant No.\ 451-03-9/2021-14/200125).

\end{document}